\documentclass[a4paper,11pt]{article}
\pdfoutput=1
\usepackage{mathrsfs}
\usepackage[utf8]{inputenc}
\usepackage[T1]{fontenc}
\usepackage{amsmath}
\usepackage{amssymb}
\usepackage{amsthm}
\usepackage[francais,english]{babel}
\usepackage{graphicx}
\usepackage{color}
\usepackage{fullpage}
\usepackage{dsfont}
\usepackage{enumerate}
\usepackage{caption}

\definecolor{purple}{rgb}{0.57,0.1,0.53}

\newenvironment{hyp}[1]{\par\bigskip \noindent \textbf{Assumption #1 : }\itshape} {\par\bigskip }
\newcommand{\HypLevy}{A}
\newcommand{\HypMut}{B.2}
\newcommand{\HypMutConst}{B.1}

\newtheorem{lemme}{Lemma}[section]
\newtheorem{theoreme}[lemme]{Theorem}

\newtheorem{proposition}[lemme]{Proposition}
\newtheorem{corollaire}[lemme]{Corollary}

\newtheorem{remarque}[lemme]{Remark}
\newenvironment{demo}{\noindent\emph{Proof :} \\ }{\hfill $\square$ \\ \par\smallskip}
\newenvironment{demoth}[1]{\noindent\emph{Proof of Theorem #1 :} \\ }{\hfill $\square$}
\newenvironment{demopr}[1]{\noindent\emph{Proof of Proposition #1 :} \\ }{\hfill $\square$}

\newcounter{claim}
\newenvironment{claim}{\stepcounter{claim}\noindent\textbf{Claim \theclaim:} \qquad }{\\}
\newenvironment*{remerciements}{
\begin{abstract}
}{\end{abstract}}

\newcommand{\Z}{\mathbb Z}
\renewcommand{\L}{\Lambda}
\newcommand{\tL}{\tilde\Lambda}

\newcommand{\sachant}{\,|\,}
\newcommand{\R}{\mathbb R}
\newcommand{\N}{\mathbb N}
\newcommand{\F}{\mathcal F}
\renewcommand{\d}[1]{\text{d}#1}
\renewcommand{\P}{\mathbb P}
\newcommand{\E}{\mathbb E}

\newcommand{\eps}{\varepsilon}

\newcommand{\D}{\mathbb D}

\newcommand{\Exc}{\mathscr E}

\newcommand{\Linv}{L^{-1}}
\newcommand{\Tx}{{T^{-x}}}
\newcommand{\Txn}{{T_{n}^{-x}}}
\renewcommand{\kill}{k}

\newcommand{\Han}{H^{+}_n}
\newcommand{\Hbn}{H^{-}_n}
\newcommand{\Hcn}{H^{\textsc{m}}_n}

\newcommand{\Ha}{H^{+}}
\newcommand{\Hb}{H^{-}}
\newcommand{\Hc}{H^{\textsc{m}}}

\newcommand{\tHa}{\tilde H^{+}}
\newcommand{\tHb}{\tilde H^{-}}
\newcommand{\tphi}{\tilde\phi}
\newcommand{\teta}{\tilde\eta}

\newcommand{\B}[1]{\mathbb B_{#1}}
\newcommand{\tr}{^t\!}

\renewcommand{\a}{\alpha}
\renewcommand{\b}{\beta}

\renewcommand{\k}{\texttt k}
\newcommand{\g}{\texttt g}
\renewcommand{\ll}{\textit{\texttt l}}
\renewcommand{\t}{\texttt t}

\newcommand{\e}{\text{e}}

\newcommand{\zero}{\{0\}}
\newcommand{\un}{\{1\}}

\newcommand{\Snk}{S_{n,k}}
\newcommand{\Snkn}{S_{n,k_n}}
\newcommand{\tLinv}{\tilde L^{-1}}

\newcommand{\Tree}{\mathbb T}
\newcommand{\tTree}{\tilde{\mathbb T}}

\newcommand{\tZ}{\tilde Z}

\newcommand{\m}{\textsc{m}}

\newcommand{\Gtaun}{G_{\tau}^{(n)}}

\newcommand{\Gtau}{G_{\tau}}

\newcommand{\gtaukn}{g_{T}^{(k_n)}}

\newcommand{\gtau}{g_T}

\begin{document}
\selectlanguage{english}
\title{Lévy processes with marked jumps I : Limit theorems}
\author{Cécile Delaporte\footnote{UPMC Univ. Paris 6, Laboratoire de probabilités et modèles aléatoires CNRS UMR 7599, 4 place Jussieu, 75252 Paris Cedex 05 France ; \textit{e-mail :} cecile.delaporte@upmc.fr ; \textit{phone :} +33144271397.}}
\date{}
\maketitle

\begin{abstract}
Consider a sequence $(\tZ_n,\tZ_n^\m)$ of bivariate Lévy processes, such that $\tZ_n$ is a spectrally positive Lévy process with finite variation, and $\tZ_n^\m$ is the counting process of marks in $\{0,1\}$ carried by the jumps of $\tZ_n$. The study of these processes is justified by their interpretation as contour processes of a sequence of splitting trees \cite{ALContour} with mutations at birth. Indeed, this paper is the first part of a work \cite{LPWM2} aiming to establish an invariance principle for the genealogies of such populations enriched with their mutational histories.

To this aim, we define a bivariate subordinator that we call the \textit{marked ladder height process} of $(\tZ_n,\tZ_n^\m)$, as a generalization of the classical ladder height process to our Lévy processes with marked jumps. Assuming that the sequence $(\tZ_n)$ converges towards a Lévy process $Z$ with infinite variation, we first prove the convergence in distribution, with two possible regimes for the marks, of the marked ladder height process of $(\tZ_n,\tZ_n^\m)$. Then we prove the joint convergence in law of $\tZ_n$ with its local time at the supremum and its marked ladder height process. The proof of this latter result is an adaptation of Chaumont and Doney \cite{CD} to the finite variation case. 
\\
\end{abstract}

\selectlanguage{french}
\section{Introduction}\label{sec_Intro}
\qquad Let $\big((Z_n,Z_n^\m)\big)_{n\geq1}$ be a sequence of bivariate Lévy processes with finite variation with values in $\R\times\Z$, such that $(Z_n,Z_n^\m)$ is characterized by its drift $(-1,0)$ and its Lévy measure $\L_n(\d r)\B{f_n(r)}(\d q)$, where $\L_n$ is a $\sigma$-finite measure on $(0,\infty)$ satisfying $\int(1\wedge r)\L_n(\d r)<\infty$, $f_n$ is a function from $(0,\infty)$ to $[0,1]$, and $\B{p}$ denotes the Bernoulli distribution with parameter $p$. We can interpret this process as a spectrally positive Lévy process with finite variation with additional marks on its jumps ; conditional on the amplitude $r$ of a jump of $Z_n$, the mark carried by this jump follows a Bernoulli distribution with parameter $f_n(r)$, and $Z_n^\m$ is then the counting process of these marks.

We consider a rescaled version $(\tZ_n,\tZ_n^\m)$ of $(Z_n,Z_n^\m)$, and assume the convergence in distribution of the sequence $(\tZ_n)$, towards a Lévy process $Z$ (with infinite variation, Assumption \HypLevy). Besides, two different assumptions concerning the marks are considered. In the first one (\HypMutConst), $(f_n)$ is a sequence of constant functions vanishing as $n\to\infty$, whereas in the second one (\HypMut), $f_n$ is a (non constant) function satisfying in particular $f_n(0)=0$. The goal of this paper is to prove some convergence theorem for the so-called \textit{marked ladder height process} of $(\tZ_n,\tZ_n^\m)$, that we define as a generalization of the classical ladder height process to Lévy processes with marked jumps. These convergence theorems are the first part of a work aiming to obtain asymptotic results for the genealogy of a splitting tree \cite{Geiger,GK,ALContour} with mutations at birth, enriched with its history of mutations. 

Let us explain how these populations can be studied from the marked Lévy processes we just described. First consider a population evolving according to the dynamics of a splitting tree $\Tree$, that is, a population where individuals give birth at constant rate during their lifetimes to i.i.d. copies of themselves. The jumping chronological contour process (or JCCP) (\cite{ALContour}) of $\Tree$ is an exploration process of this tree that provides a one-to-one correspondence with $\Tree$, and which distribution is characterized from a spectrally positive Lévy process with finite variation. Assume now that individuals carry types, and that (neutral) mutations may happen at birth of individuals : to each birth event in $\Tree$ we associate a mark in $\{0,1\}$, which will code for the absence ($0$) or presence ($1$) of a mutation. Then the generalization of the JCCP for this splitting tree with marks leads to a characterization of its law by a spectrally positive Lévy process with finite variation, with 
additional 
marks on its jumps, as described earlier.\\

Thus let us interpret our sequence $\big((Z_n,Z_n^\m)\big)$ as the contour processes of a sequence of marked splitting trees $(\Tree_n)$. Roughly speaking, the measure $\L_n$ characterizes the lifetime distribution of the individuals in $\Tree_n$, and conditional on its lifetime $r$, an individual has probability $f_n(r)$ to be a mutant. Letting $n\to\infty$, we aim at stating results in a large population asymptotic for $\Tree_n$, which requires to introduce a rescaling of these populations. Here the convergence assumption on $(\tZ_n)$ has to be interpreted as the convergence, in a certain sense, of the populations $(\tTree_n)$ obtained from a proper rescaling of $(\Tree_n)$. 

More precisely, our ultimate goal is to obtain an invariance principle for the genealogy (with mutational history) of the rescaled population $\tTree_n$, as $n\to\infty$. The characterization of the latter with the help of the JCCP can be obtained from the law of the (marked) future infimum of an excursion of the Lévy process $\tZ_n$ under a fixed level. By a time reversal argument, this comes to study the (marked) running supremum of $\tZ_n$ killed upon hitting $0$. Here \og marking\fg\ the future infimum (resp. running supremum) of $\tZ_n$ means selecting and keeping record of the marks carried by the jumps of the future infimum (resp. running supremum) of $\tZ_n$. We are thus led to introduce the marked ladder height process of $(\tZ_n,\tZ_n^\m)$ : consider $\Ha_n$ the ascending ladder height process of $\tZ_n$, and put marks on its jumps in agreement with the marks on the corresponding jumps of $\tZ_n$. Denoting by $\Hc_n$ the counting process of these marks, the so-called marked ladder height process $(\
Ha_n,\Hc_n)$ is 
then a (possibly killed) bivariate subordinator. 

We are here interested in the asymptotic behaviour of these processes under Assumptions \HypLevy\ and \HypMutConst/\HypMut\ defined above (see also Section \ref{sec_cv_ass}). While Assumption \HypLevy\ alone ensures the convergence in distribution of $\Ha_n$ towards the classical ladder height process of $Z$, Assumptions \HypMutConst\ and \HypMut\ are designed to allow that of the marked ladder height process. We prove in Section \ref{sec_cvMLHP} the convergence in law of $(\Ha_n,\Hcn)$ towards a (possibly killed) bivariate subordinator $(\Ha,\Hc)$, such that $\Ha$ is the ladder height process of $Z$. Note nevertheless that in this framework there is in general no convergence of the whole mutation process, namely $\tZ_n^\m$. 
In the case of Assumption \HypMutConst, $\Ha$ and $\Hc$ are independent, and $\Hc$ is a Poisson process with parameter $\theta$, which arises as the limit of the sequence of constant functions $(f_n)$ after a proper rescaling. This means that the contribution to the marks in the limit exclusively comes from jumps with vanishing amplitudes. This is no longer the case under Assumption \HypMut, yet additional independent marks can appear if $Z$ has a Gaussian component. In Section \ref{sec_jointcv} we establish the joint convergence in law of $(\tZ_n,L_n,\Ha_n,\Hc_n)$, where $L_n$ is a local time of $\tZ_n$ at its supremum. The proof of this result is essentially an adaptation of L. Chaumont and R.A. Doney's paper \cite{CD}, to our specific case of finite variation Lévy processes converging to an infinite variation Lévy process. 

\setlength{\parindent}{0cm}

\section{Preliminaries}\label{sec_Prelim}
This section sets up notation for the topological framework, and provides some background on spectrally positive Lévy processes.
\subsection{Topology} \label{sec_topo} 

We consider the Euclidean space $\R^d$ and endow it with its Borel $\sigma$-field $\mathcal B(\R^d)$. For all $x\in\R^d$, $\tr x$ will denote its transpose. We denote by $\D(\R^d)$ the space of all càd-làg functions from $\R_+$ to $\R^d$. We endow the latter with the Skorokhod topology, which makes it a Polish space (see \cite[VI.1.b]{JS}). In the sequel, for any function $f\in\D(\R)$ and $x>0$, we will use the notation $\Delta f(x)=f(x)-f(x-)$, where $f(x-)=\lim_{u\to x,\,u<x} f(u)$.\\

Now for any Polish space $X$, with its Borel $\sigma$-field $\mathcal B$, the space $\mathcal M_f(X)$ of positive finite measures on $(X,\mathcal B)$ can be endowed with the weak topology :
It is the coarsest topology for which the mappings $\mu\mapsto\int g\d\mu$ are continuous for any continuous bounded function $g$. In the sequel, we will use the notation $\mu(g):=\int g\d\mu$.\\

Hence we endow here $\mathcal M_f(\R^d)$ and $\mathcal M_f(\D(\R^d))$ with their respective weak topologies. The notation $\Rightarrow$ will be used for both weak convergence in $\R^d$ and in $\D(\R^d)$, and we will use the symbol $\stackrel{(d)}{=}$ for the equality in distribution. Recall that for any sequence of $\R^d$-valued càd-làg processes $(X_n)$, the weak convergence of $(X_n)$ towards a process $X$ of $\D(\R^d)$ is equivalent to the finite dimensional convergence of $(X_n)$ towards $X$ along any dense subset $D\subset\R_+$, together with the tightness of $(X_n)$. For more details about convergence in distribution in $\D(\R^d)$, see \cite[VI.3]{JS}.\\

\subsection{Spectrally positive Lévy processes} \label{sec_SPLP}
This paragraph is composed of results that can mostly be found in \cite{B} or \cite{K}, and consists in a summary of the main points concerning spectrally positive Lévy processes.\\

We consider a real-valued Lévy process $X=(X_t)_{t\geq0}$ (that is, $X$ is a càd-làg process with independent and stationary increments), which we will suppose spectrally positive, meaning that it has no negative jumps. We assume furthermore that $X$ is starting at $0$ a.s., and denote by $\P$ its law. This Lévy process is characterized by its Laplace exponent $\psi$ defined for all $\lambda\geq0$ by
$$\E(e^{-\lambda X_t})=e^{t\psi(\lambda)},$$
and the Lévy-Khintchine formula gives :
\begin{equation} \label{formula_LK}
 \psi(\lambda):=d\lambda+\frac{b^2}2\lambda^2-\int_{(0,\infty)}(1-e^{-\lambda r}-\lambda h(r))\L(\d r),
\end{equation}
where $h$ is some arbitrary truncation function on $\R$ (in general, a truncation function $h$ is a continuous bounded function from $\R^d$ to $\R^d$ satisfying $h(x)\sim x$ in a neighbourhood of 0). The Lévy measure $\L$ is a measure on $(\R_+^*,\mathcal B(\R_+^*))$ satisfying $\int (1\wedge |u|^2) \L(\d u)<\infty$. The coefficient $b$ is named Gaussian coefficient, and the coefficient $d$ depends on the choice of the truncation function.\\

The paths of $X$ have finite variation (on every compact time interval) a.s. iff $b=0$ and $\int (1\wedge |r|)\L(\d r)<\infty$. In this case, the integral $\int_{(0,\infty)} h(r)\L(\d r)$ is finite a.s., and we can reexpress the Laplace exponent as
\begin{equation} \label{formula_LK_vf}
 \psi(\lambda):=-d'\lambda-\int_{(0,\infty)}(1-e^{-\lambda r})\L(\d r),
\end{equation}
where $d'$ is called the drift coefficient and characterizes $X$ together with the Lévy measure $\L$. It is in particular the case if $X$ is a subordinator, i.e. if $X$ has increasing paths a.s., and then $d'$ is nonnegative. In the sequel, we will sometimes deal with killed subordinators : by killed subordinator at a random time $T$ we mean that the value of the process at any time $t\geq T$ is replaced by $+\infty$. By killed subordinator at rate $k$ we mean a killed subordinator at an independent exponential time with parameter $k$. \\

Consider the case where $X$ is not a subordinator (note that if $X$ has finite variation, it has necessarily a drift $d'<0$). The Laplace exponent $\psi$ is infinitely differentiable, strictly convex, and satisfies $\psi(0)=0$ and $\underset{\lambda\to\infty}\lim \psi(\lambda)=+\infty$. In particular, $\psi'(0^+)=-\E(X_1)\in[-\infty,+\infty)$. Thus $\psi$ has at most one root besides $0$. We denote by $\eta$ the largest one, and $\eta=0$ if and only if $\psi'(0^+)\geq0$. Moreover, $X$ drifts to $+\infty$ (resp. oscillates, drifts to $-\infty$) if and only if $\psi'(0^+)$ is negative (resp. zero, positive). Then we say that $X$ is respectively supercritical, critical or subcritical. Note that if $X$ is supercritical, $\eta>0$, and that otherwise $\eta=0$. Furthermore, the function $\psi$ is a bijection from $[\eta,\infty)$ to $\R_+$ and we define its inverse $\phi:\ \R_+\to[\eta,\infty)$. \\

Finally we introduce the scale function, which is in particular useful for solving exit problems (see e.g. \cite[Chapter 8]{K}) : $W$ is defined as the unique strictly increasing continuous function from $\R_+$ to $\R_+$ with Laplace transform
\begin{equation}\label{formula_W}
 \int_{(0,\infty)} e^{-\lambda x} W(x) \d x=\frac{1}{\psi(\lambda)},\ \ \lambda>\eta.
\end{equation}

According to \cite[Lemma 8.6]{K}, when $X$ is not a subordinator, $W(0)$ is equal to $-1/d'$ (where $d'<0$ is the drift) in the finite variation case, and is zero in the infinite variation case.

\subsection{Local time and excursions}\label{sec_Exc}
Let $X$ be a spectrally positive Lévy process with Laplace exponent $\psi$ given by formula (\ref{formula_LK}), and denote by $(\F_t)$ the natural filtration associated with $X$, i.e. for all $t\geq0$,
$$\F_t=\sigma\{X_s,\ s\leq t\}.$$

We define its past supremum $\bar X_t:=\underset{[0,t]}\sup X$ for all $t\geq0$. Then the reflected process $X-\bar X$ is a Markov process in the filtration $(\F_t)$ (and also in its own natural filtration), for which one can construct a local time at $0$ and develop an excursion theory. For more details about the following results, see chapter IV in \cite{B}.\\

\paragraph*{Local times} \quad

For the construction of a local time at $0$ for $X-\bar X$ (which we will also name local time at the supremum for $X$), we have to distinguish the case of infinite variation, where $0$ is regular for $X$ w.r.t. the open half-line $(0,\infty)$, from the case of finite variation, where $0$ is irregular w.r.t. the open half-line $(0,\infty)$. 

According to Theorem IV.4 in \cite{B}, when $X$ has infinite variation, we denote by $L$ a local time at $0$ for $X-\bar X$, and the mapping $t\mapsto L(t)$ is non decreasing and continuous. Any other local time at $0$ for $X-\bar X$ differs then from $L$ in a positive multiplicative constant. When $X$ has finite variation, we set 
\[L(t):=\sum_{i=0}^{\ll(t)} \tau_i,\]
where $\ll(t)$ represents the number of jumps of the supremum up until time $t$ - i.e. the number of zeros of the reflected process up until time $t$, and $(\tau_i)_{i\geq0}$ is a sequence of i.i.d. random exponential  variables with arbitrary parameter, independent from $X$. Then $L$ is a local time at the supremum for $X$, but is only right-continuous. However, $L$ is not adapted to the filtration $(\F_t)$, and to make up for that problem we replace $(\F_t)$ by $(\mathcal G_t):=(\F_t \vee \sigma(L_s,\,s\leq t))$. We can then define in both cases the right-continuous inverse of $L$ : for all $t\geq0$, set 
$$\Linv(t):=\inf\{s\geq0,\ L(s)>t\}.$$

The process $\Linv$ is a killed subordinator, and is adapted to $(\mathcal G_{\Linv(t)})$.

\paragraph*{Excursion theory} \quad

We denote by $\Exc$ the set of excursions of $X-\bar X$ away from $0$ : $\Exc$ is the set of the càd-làg functions $\epsilon$ with no negative jumps for which there exists $\zeta=\zeta(\epsilon)\in(0,\infty]$, which will be called the lifetime of the excursion, and such that $\epsilon(0)=0$, $\epsilon(t)$ has values in $(-\infty,0)$ for $t\in(0,\zeta)$ and in the case where $\zeta<\infty$, $\epsilon(\zeta)\in[0,\infty)$.\\

We consider the process $e=(e_t)_{t\geq0}$ with values in $\Exc\cup \{\partial\}$ (where $\partial$ is an additional isolated point), defined by :
\[e_t:=\left\{ 
\begin{array}{l l}
  ((X-\bar X)_{s+\Linv(t-)},0\leq s<\Linv(t)-\Linv(t-)) & \quad \text{if } \Linv(t-)<\Linv(t)\\
  \partial & \quad \text{else}\\ \end{array} \right. .\]

Then according to Theorem IV.10 in \cite{B}, if $X$ does not drift to $-\infty$, then $0$ is recurrent for the reflected process, and $(t,e_t)_{t\geq0}$ is a Poisson point process with intensity $c\;\d t\ N(\d\epsilon)$, where $c$ is some constant depending on the choice of $L$, and $N$ is a measure on $\Exc$. Else, $(t,e_t)_{t\geq0}$ is a Poisson point process with intensity $c\;\d t\ N(\d\epsilon)$, stopped at the first excursion with infinite lifetime. \\

Finally, we describe some marginals of $N$ in the proposition below, for which we refer to \cite[Th. 6.15 and (8.29)]{K}, \cite[(3)]{BertoinDecomp} and \cite[Cor. 1]{BertoinPitman}.

\begin{proposition} \label{prop_Exc_mesure} 
 We have for all $z,x>0$ :
\begin{enumerate}[\upshape(i)]
 \item If $X$ has finite variation,
\upshape$$N(-\epsilon(\zeta-)\in \d x,\ \epsilon(\zeta)\in \d z,\ \zeta<\infty)=W(0)e^{-\eta x}\d x\L(x+\d z).$$
\itshape
 \item  If $X$ has infinite variation and no Gaussian component (i.e. $b=0$),
\upshape$$N(-\epsilon(\zeta-)\in \d x,\ \epsilon(\zeta)\in \d z,\ \zeta<\infty)=e^{-\eta x}\d x\L(x+\d z).$$
\end{enumerate}
\itshape Moreover, in both cases, under $N(\,\cdot\,\sachant-\epsilon(\zeta-)=x,\ \zeta<\infty)$, the reversed excursion 
$$\big(-\epsilon((\zeta-t)-),\ 0\leq t<\zeta\big)$$ 
is equal in law to $(X_t,\ 0\leq t<T^0)$ under $\P_x(\,\cdot\,\sachant T^0<\infty)$.
\end{proposition}

\paragraph*{Convergence of Lévy processes}\quad

Finally, we recall a restricted version of Corollary 3.6 from \cite[VII.3]{JS}, that will be needed later : 

\begin{proposition}\label{prop_thJS}
Let $X_n$, $X$ be spectrally positive Lévy processes with respective Laplace exponents
\upshape
$$\psi_n(\lambda):=c_n\lambda+\frac{b_n^2}2 \lambda^2-\int (1-e^{-\lambda u}-\lambda h(u))\L_n(\d u)$$
$$\psi(\lambda):=c\lambda+\frac{b^2}2 \lambda^2-\int (1-e^{-\lambda u}-\lambda h(u))\L(\d u)$$
\itshape
for some common truncation function $h$. Then $X_n\Rightarrow X$ in $\D(\R_+)$ iff as $n\to\infty$ :
\begin{enumerate}[\upshape (i)]
 \item $c_n\to c$, \label{thJS_drift}
 \item \upshape  $b_n^2+\int h^2 \d\L_n \to b^2+\int h^2 \d\L$, \itshape\label{thJS_coeff_brownien}
 \item \itshape For any continuous bounded function $g$ satisfying $g(u)=o(|u|^2)$ when $|u|\to0$ (or equivalently, vanishing on a neighbourhood of 0),\upshape $\int g\d\L_n \to \int g \d\L$. \label{thJS_mesure_Levy}
\end{enumerate}
\end{proposition}

\begin{remarque}
An analogous version of this statement is available for Lévy processes with values in $\R^d$, for which each coordinate is itself spectrally positive. Note in particular that condition \eqref{thJS_coeff_brownien} is then : for all $1\leq i,j\leq d$, \upshape $b_n^2+\int h_i h_j \d\L_n \to b^2+\int h_i h_j \d\L$\itshape, where $h_i$ denotes the $i$-th coordinate of $h$.
\end{remarque}

\subsection{Lévy process with marked jumps and marked ladder height process}\label{sec_mLHP_def}

Let $\L$ be a measure on $(\R_+^*,\mathcal B(\R_+^*))$ satisfying $\int(1\wedge u) \L(\d u)<\infty$, and $f$ a function from $\R_+^*$ to $[0,1]$. Denote by $\B r$ the Bernoulli probability measure with parameter $r$, and consider $(X,X^\m)$ a bivariate Lévy process with finite variation, with Lévy measure $\L(\d u)\B{f(u)}(\d q)$ and drift $(-1,0)$. These marked Lévy processes will be used in \cite{LPWM2} to characterize the law of the contour of a splitting tree with mutations at birth, as explained in Section \ref{sec_Intro}. We define now the marked ladder height process of $X$. This process is a bivariate subordinator, whose first coordinate will be the classical ladder height process of $X$, and whose second coordinate will keep record of the marks that are present on the jumps of the current supremum of $X$. It appears naturally in the second paper \cite{LPWM2}, as a tool to describe the distribution of mutations on the genealogy of a marked splitting tree. \\

Sticking to the notation introduced in Section \ref{sec_SPLP} for $X$ and in Section \ref{sec_Exc} for the local time and excursion process of $X-\bar X$, we define for all $t\in[0,L(\infty))$
 $$\xi_t:=\left\{ 
\begin{array}{l l}
 (t,e_{t}(\zeta),-e_{t}(\zeta-),\Delta X^\m(\Linv(t))) & \quad \text{if } \Linv(t-)<\Linv(t)\\
 \partial & \quad \text{else}\\
\end{array}
\right.,$$
where $\partial$ is an additional isolated point, and $e_{t}(\zeta)$ (resp. $e_{t}(\zeta-)$) stands for $e_{t}(\zeta(e_{t}))$ (resp. $e_{t}(\zeta(e_{t})-)$).

Here the fourth coordinate $\Delta X^\m(\Linv(t))$ is $1$ or $0$ whether or not the jump of $X$ at the right end point of the excursion interval indexed by $t$ carries a mark. Note that the set $\{\Linv(t)\}_{t\geq0}$ of these right end points is exactly the set of record times of $X$.

\begin{lemme} \label{lemme_xi}
The process $\xi$ is distributed as a Poisson point process on $[0,\mathcal K)\times\R_+^*\times\R_+^*\times\{0,1\}$ with intensity measure
\upshape$$cW(0)\,\d t \cdot  \L(x+\d y)\; e^{-\eta x} \d x   \cdot \B{f(x+y)}(\d q),$$
\itshape where if $X$ drifts to $-\infty$, $\mathcal K$ is an independent exponential variable with parameter $\kill:=c\frac{W(0)}{W(\infty)}$, and else $\mathcal K=+\infty$ a.s. 
\end{lemme}

\begin{demo}
We denote by $\tilde\xi$ the restriction of $\xi$ to its first three coordinates. \\
We know from \cite[Prop. 0.5.2]{B} and Section \ref{sec_Exc} that $\tilde\xi$ is distributed as a Poisson point process on $[0,\mathcal K)\times\R_+^*\times\R_+^*$ with intensity 
$$c \d t \  N(\epsilon(\tilde\xi)\in\d y,\,-\epsilon(\tilde\xi-)\in\d x),$$ 
where from Proposition \ref{prop_Exc_mesure},
$$N(\epsilon(\tilde\xi)\in\d y,\,-\epsilon(\tilde\xi-)\in\d x)=W(0) e^{-\eta x}\d x\L(x+\d y),$$
and $\mathcal K$ is an independent exponential variable with parameter $c N(\{\epsilon\in\Exc,\ \tilde\xi(\epsilon)=\infty\})=c\frac{W(0)}{W(\infty)}$ if $X$ drifts to $-\infty$, and else $\mathcal K=+\infty$ a.s. \\
Let $B\in\mathcal B(\R^*_+\times\R_+^*)$, and $t\geq0$. Conditional on having an atom of $\tilde\xi$ in $[0,t]\times B$, the fourth coordinate of the corresponding atom of $\xi$ follows a Bernoulli distribution with parameter :
\begin{align*}
 p(B):=\frac{\int_B f(x+y)\ N(\d y,\d x)}{N(B)}.
\end{align*}
As a consequence, $\xi([0,t]\times B\times\un)$ and $\xi([0,t]\times B\times\zero)$ follow Poisson distributions with respective parameters $p(B)N(B)c t$ and $(1-p(B))N(B)c t$, and we deduce that $\xi$ is a Poisson random measure with intensity $\pi$, such that for $C\in\mathcal P(\{0,1\})$ :
\begin{align*}
\pi([0,t]&\times B\times C)\\
= & c t\ N(B)\B{p(B)}(C) \\
= & c t\ \int_B \B{f(x+y)}(C) N(\d y,\d x),
\end{align*}
which leads to the result.
\end{demo} 

Let $(\Ha,\Hb,\Hc)$ be the (possibly killed) trivariate subordinator with no drift and whose jump point process is a.s. equal to the restriction of $\xi$ to its last three coordinates.  Here we define $\Hb$ only for technical reasons (see Section  \ref{sec_jointcv}), and hence we now define the marked ladder height process of $X$ as the (possibly killed) bivariate subordinator $(\Ha,\Hc)$. However it will be convenient in the sequel to be also able to name $(\Ha,\Hb,\Hc)$ ; we call it the trivariate ladder height process of $X$.\\

Then, as a straightforward consequence of Lemma \ref{lemme_xi} we have 
\begin{proposition}\label{prop_lhp}
 The marked ladder height process  $(\Ha,\Hc)$ is a bivariate subordinator with no drift and Lévy measure
\upshape\begin{equation} cW(0) \int_0^\infty  \d x\; e^{-\eta x}\; \L(x+\d y) \; \B{f(x+y)}(\d q),
\end{equation}\itshape 
and killed at rate $\kill=c\frac{W(0)}{W(\infty)}$.
\end{proposition}

 Note that $\Ha$ is in fact the ladder height process of $X$, i.e. for all $t\geq0$, $\Ha(t)=\bar X(\Linv(t))$ a.s. Moreover, $\Hc$ is a Poisson process which jumps correspond, in the local time scale, to the marks occurring at record times of $X$.

\section{Definitions and notation}

\subsection{Convergence assumptions}\label{sec_cv_ass}

Let $(\L_n)_{n\geq1}$ be a sequence of measures on $(\R_+^*,\mathcal B(\R_+^*))$ satisfying $\int(1\wedge u) \L_n(\d u)<\infty$ for all $n$, and $(f_n)_{n\geq1}$ a sequence of continuous functions from $\R^+$ to $[0,1]$. We consider a sequence of independent bivariate Lévy processes $(Z_n,Z_n^{\textsc m})_{n\geq1}$ with finite variation, Lévy measure $\L_n(\d u)\B{f_n(u)}(\d q)$ and drift $(-1,0)$, where we recall that $\B r$ denotes the Bernoulli probability measure with parameter $r$. We first assume 

 \begin{hyp}\HypLevy
  There exists a sequence of positive real numbers $(d_n)_{n\geq1}$ such that as $n\to\infty$, the process defined by
$$\tilde Z_n:=\Big(\frac1n Z_{n}(d_n t)\Big)_{t\geq0}$$ 
converges in distribution to a (necessarily spectrally positive) Lévy process $Z$ with infinite variation, and with Lévy measure denoted by $\L$.
 \end{hyp}
For all $n\in\N$ and $t\geq0$, set $\tilde Z_n^{\textsc m}(t):=Z_n^{\textsc m}(d_n t)$. In the sequel we always assume that $\tZ_n(0)=\tZ_n^\m(0)=0$. With a slight abuse of notation, the law of $(\tilde Z_n,\tilde Z_n^{\textsc m})$ conditional on $(\tilde Z_n(0),\tilde Z_n^\m(0))=(0,0)$, and the law of $Z$ conditional on $Z(0)=0$, will both be denoted by $\P$.

\paragraph*{Some notation :} 
As in Section \ref{sec_SPLP}, the Laplace exponents $\psi_n$ of $Z_n$, $\tilde \psi_n$ of $\tilde Z_n$ and $\psi$ of $Z$ are defined by 
$$\E(e^{-\lambda Z_n(t)})=e^{t\psi_n(\lambda)},\  \E(e^{-\lambda \tilde Z_n(t)})=e^{t\tilde\psi_n(\lambda)}\  \text{ and }\ \E(e^{-\lambda Z(t)})=e^{t\psi(\lambda)},\ \ \ \lambda\geq0.$$
We denote by $\teta_n$ (resp. $\eta$) the largest root of $\tilde\psi_n$ (resp. $\psi$) and by $\tilde\phi_n$ (resp. $\phi$) the inverse of $\tilde\psi_n$ (resp. $\psi$) on $[\teta_n,\infty)$ (resp. $[\eta,\infty)$). We denote by $\tilde W_n$ (resp. $W$) the scale function of $\tilde Z_n$ (resp. $Z$). Finally, we denote by $\tL_n$ the Lévy measure of $\tilde Z_n$.

\paragraph*{Remarks about $(d_n)$ :} 
Writing for $\lambda\geq0$, $\E(e^{-\lambda \tilde Z_n(t)})=e^{d_n t\psi_n(\lambda/n)}$, we get from formula (\ref{formula_LK_vf})  that $\tilde Z_n$ has drift $-\frac{d_n}n$, Lévy measure $\tL_n=d_n\L_n(n\cdot)$ and Laplace exponent $\tilde\psi_n=d_n\psi(\cdot/n)$. In particular, this gives $\tilde W_n(0)=n/d_n$. We state later in Proposition \ref{prop_cv} that $\tilde W_n$ converges pointwise to $W$ as $n\to\infty$, and besides, the assumption of infinite variation of $Z$ ensures $W(0)=0$. Thereby we know that necessarily $\frac{d_n}{n}\to\infty$ as $n\to\infty$.\\

\par\bigskip
Finally, we suggest two possible assumptions for the asymptotic of the marks : in the first one, the probability for a jump of $\tZ_n$ to carry a mark is constant, while in the second one, this probability is a function of the amplitude of the jump.

\begin{hyp}\HypMutConst
\begin{enumerate}[\upshape(a)]
 \item For all $n\geq1$, for all $u\in\R_+$, $f_n(u)=\theta_n$, where $\theta_n\in[0,1]$.
 \item As $n\to\infty$, $\frac{d_n}{n}\theta_n$ converges to some finite real number $\theta$.
 \end{enumerate}
\end{hyp}

\begin{hyp}\HypMut
\begin{enumerate}[\upshape(a)]
 \item The sequence $\big(u\mapsto \frac{f_n(nu)}{1\wedge u}\big)$ converges uniformly to $u\mapsto \frac{f(u)}{1\wedge u}$ on $\R_+^*$. \label{HypMut_cv}
 \item There exists $\kappa\geq0$ such that $f(u)/u \to \kappa$ as $u\to0^+$. \label{HypMut_0} 
\end{enumerate}
\end{hyp}

Note that in \HypMutConst, necessarily $\theta_n\to 0$ as $n\to\infty$. Then if we denote by $f$ the limit of the sequence $(f_n)$, we have $f\equiv0$. Besides, in Assumption \HypMut\ the choice of $f_n$ and $f$ is independent of $\tilde Z_n$ and $Z$. 

\begin{remarque}\label{remark_pas_de_cv_Zm}
These two possible assumptions have been chosen so that as $n\to\infty$, we have convergence of the set of marks that are carried by jumps of the supremum (which will be reexpressed as sets of mutations on a lineage in the second paper \cite{LPWM2}). However this choice does not imply, despite Assumption \HypLevy, the convergence of the bivariate process $(\tilde Z_n,\tilde Z_n^\m)$. It is even never the case under \HypMut\ : from Proposition \ref{prop_thJS} we see that the convergence as $n\to\infty$ of \upshape$\int_{(0,\infty)} f_n(nu)\tL_n(\d u)$\itshape is a necessary condition for that of $(\tilde Z_n,\tilde Z_n^\m)$. Now it can be shown that under \HypMut, this integral behaves as $n\to\infty$ like \upshape$\int_{(0,\infty)} (1\wedge u)\tL_n(\d u)$\itshape, which goes to $\infty$ as $n\to\infty$ (see Lemma \ref{lemme_theta_n} for a similar result).
\end{remarque}

\subsection{Marked ladder height process of $\tZ_n$} 

\paragraph*{Local times at the supremum}\quad 

We denote by $\F=(\F_t)_{t\geq0}$ (resp. $\F_n=(\F_{n,t})_{t\geq0}$) the natural filtration associated to $Z$ (resp. $\tZ_n$), that is for all $t\geq0$,
$$\F_t=\sigma\{Z_s,\ s\leq t\}\ (\text{resp. }\ \F_{n,t}=\sigma\{\tZ_n(s),\ s\leq t\}).$$

For all $n\geq1$, let $(\tau_{n,i})_{i\geq0}$ be a sequence of i.i.d. random exponential variables, independent of $(\tZ_n)_{n\geq1}$, with parameter $\alpha_n:=\frac{d_n}{n}$. This choice will allow us in the sequel to obtain some convergence properties, in particular for the inverse local time and the ladder height process of $\tZ_n$. Then, according to Section \ref{sec_SPLP}, we define for $\tZ_n$ a local time at the supremum as follows :
\[L_n(t):=\sum_{i=0}^{\ll_n(t)} \tau_{n,i},\]
where $\ll_n(t)$ represents the number of jumps of the supremum until time $t$. We denote by $\Linv_n$ the right-continuous inverse of $L_n$ as defined in Section \ref{sec_SPLP}, and replace the filtration $\F_{n,t}$ with $\F_{n,t} \vee \sigma (L_n(s),\, s\leq t)$, so that $L_n$ (resp. $\Linv_n$) is adapted to $(\F_{n,t})$ (resp. to $(\F_{n,\Linv_n(t)})$).\\

As in Section \ref{sec_Exc}, we introduce the local time at the supremum $L$ for the infinite variation Lévy process $Z$ : we saw that $L$ is defined up to a multiplicative constant, and we require that
\begin{equation} \label{formula_norm_local_time}
\E\bigg(\int_{(0,\infty)}e^{-t}\d L_t\bigg)=\phi(1),
\end{equation}
so that $L$ is uniquely determined. Finally, we denote by $\Linv$ its inverse.\\

\paragraph*{Marked ladder height process} \quad

For $n\geq1$, let $(\Han,\Hbn,\Hcn)$ be the trivariate marked ladder height process of $\tZ_n$, as defined in section \ref{sec_mLHP_def}. Recall that we are mostly interested in $(\Han,\Hcn)$ and that we define $\Hbn$ only for technical reasons (see Section \ref{sec_jointcv}). For this reason in the sequel, we focus on $(\Han,\Hcn)$. The results will first be stated in terms of the (bivariate) ladder height process, but their proofs can be easily adapted to the trivariate ladder height process.

Our choice for the normalization of the local times, and the equality $\tilde W_n(0)=\frac{n}{d_n}$, along with Proposition \ref{prop_lhp}, yields
\begin{proposition}
 The ladder height process $(\Han,\Hcn)$ is a bivariate subordinator with no drift and Lévy measure \upshape
\begin{equation} \label{mu_n}
\mu_n(\d y,\d q):=\int_0^\infty  \d x\; e^{-\teta_n x}\; \tL_n(x+\d y) \; \B{f_n(n(x+y))}(\d q),
\end{equation} \itshape
and killed at rate $\kill_n:=\frac{1}{\tilde W_n(\infty)}$ if $\tZ_n$ is subcritical.
\end{proposition}

We also introduce the notation 
\begin{equation} \label{mu+_n}
\mu^+_n(\d y):=\mu_n(\d y,\{0,1\})=\int_0^\infty  \d x\; e^{-\teta_n x}\; \tL_n(x+\d y)
\end{equation}
for the Lévy measure of $\Han$. As stated in Section \ref{sec_mLHP_def}, $\Han$ is in fact the ladder height process of $\tZ_n$, i.e. for all $t\geq0$, $\Han(t)=\bar \tZ_n(\Linv_n(t))$ a.s., where $\bar \tZ_n(t)$ denotes the current supremum of $\tZ_n$ at time $t$. Moreover, $\Hcn$ is a Poisson process with parameter $\lambda_n:=\mu_n(\R_+^*\times\un)$, so that the random time 
\begin{equation} \label{e_n}
\e_n:=\inf\{t\geq0,\ \Hcn(t)=1\}
\end{equation}
follows on $\{\e_n<L_n(\infty)\}$ an exponential distribution with parameter $\lambda_n$.

\section{Convergence theorem for the marked ladder height process}\label{sec_cvMLHP}

\subsection{Statement of result}
We define 
$$\mu(\d u,\d q):=\int_0^\infty  \d x\; e^{-\eta x}\; \L(x+\d u) \; \B{f(x+u)}(\d q),$$
and
$$\mu^+(\d u):=\mu(\d u,\{0,1\})=\int_0^\infty  \d x\; e^{-\eta x}\; \L(x+\d u).$$

Then, we have the following theorem :
\begin{theoreme} \label{th_cv_H}
Under Assumption \HypMutConst, if $Z$ does not drift to $-\infty$, the sequence of bivariate subordinators $H_n=(\Han,\Hcn)$ converges weakly in law to a subordinator $H:=(\Ha,\Hc)$, where $\Ha$ and $\Hc$ are independent, $\Ha$ is a subordinator with drift $\frac{b^2}2$ and Lévy measure $\mu^+$, and $\Hc$ is a Poisson process with parameter $\theta$.
In the case $Z$ drifts to $-\infty$, the same statement holds but $H$ is killed at rate $\kill:=\frac1{W(\infty)}$ and the independence between $\Ha$ and $\Hc$ holds only conditional on their common lifetime.
\\

 Under Assumption \HypMut, the sequence of bivariate subordinators $H_n=(\Han,\Hcn)$ converges weakly in law to a subordinator $H:=(\Ha,\Hc)$, which is killed at rate $\kill$ if $Z$ drifts to $-\infty$. Moreover, $H$ has drift $(\frac{b^2}2,0)$ and Lévy measure 
\upshape $$\mu(\d u,\d q)+\rho\delta_0(\d u)\delta_1(\d q),$$
\itshape where $\rho:=\kappa b^2$. 
\end{theoreme}
In particular, under Assumption \HypMut, if $Z$ has no Gaussian component, the limiting marked ladder height process is a pure jump bivariate subordinator with Lévy measure $\mu$. If $Z$ has a Gaussian component, the fact that the \og small jumps\fg\ of $\tZ_n$ generate the Gaussian part in the limit results in a drift for $\Ha$, and possibly additional independent marks that happen with constant rate in time, as under Assumption \HypMutConst. This rate is proportional to the Gaussian coefficient (provided that $\kappa\neq0$). Besides, note that as expected, $\Ha$ is distributed as the classical ladder height process of $Z$. The joint convergence in law of the triplet $(\tZ_n,\Han,\Hcn)$ towards $(Z,\Ha,\Hc)$ is established in the next section.
\par\bigskip
\begin{remarque} \label{remark_cv_H-}
For technical reasons we also need to obtain the convergence in distribution of $(\Han,\Hbn,\Hcn)$. According to Lemma \ref{lemme_xi}, this process is a trivariate pure jump subordinator with Lévy measure 
\upshape$$  \d x\; e^{-\teta_n x}\; \tL_n(x+\d u) \; \B{f_n(n(x+u))}(\d q),$$\itshape
and we can easily adapt the upcoming proofs to get that $(\Han,\Hbn,\Hcn)$ converges in distribution to a subordinator $(\Ha,\Hb,\Hc)$.
\end{remarque}

\subsection{Proof}
\paragraph*{Consequences of Assumption \HypLevy}\quad

Before proving Theorem \ref{th_cv_H}, we state some direct consequences of the convergence of $\tZ_n$ towards $Z$. The two following propositions will be frequently used in the sequel and shall be kept in mind by the reader.\\

\begin{proposition}\label{prop_cv}
\begin{enumerate}[\upshape(i)]
 \item As $n\to\infty$, $\tilde\phi_n\to\phi$ uniformly on every compact set of $\R_+$, and in particular $\tilde\eta_n\to\eta$. \label{prop_cv_phi}
 \item As $n\to\infty$, $\tilde W_n\to W$ uniformly on $\R_+$. \label{prop_cv_W}
\end{enumerate}
\end{proposition}

\begin{demo}
Denote by $T_n^x$ (resp. $T^x$) the first entrance time of $\tZ_n$ (resp. $Z$) in the Borel set $\{x\}$, $x\in\R$. Since $Z$ has no negative jumps it is a.s. continuous at $\Tx$, and we have $\underset{\eps\to0+}\lim T^{-(x+\eps)}=\Tx$ a.s. Hence as a straightforward consequence of Proposition VI.2.11 in \cite{JS}, we have the convergence in law of $\Txn$ towards $\Tx$. 
Now $\phi_n$ (resp. $\phi$) is the Laplace exponent of the process $x\mapsto \Txn$ (resp. $x\mapsto \Tx$) \cite[Th. VII.1.1]{B}. The pointwise convergence of $\tilde\phi_n$ to $\phi$ is thus a consequence of the convergence in distribution of $\Txn$ towards $\Tx$. The uniform convergence comes from the fact that for all $n\geq1$, $\tilde\phi_n$ is increasing on $\R_+$.\\
The proof of the pointwise convergence of $\tilde W_n$ towards $W$ can be found in \cite[Prop. 3.1]{LS} or can be derived from its definition. Moreover, we have for all $y>x$ $\P(T^{-x}<T^{(y-x,\infty)})=\frac{W(x)}{W(y)}$ \cite[Th. VII.2.8]{B}, and then the function $x\mapsto \tilde W_n(x)/\tilde W_n(y)$ is decreasing. The convergence of $\tilde W_n$ towards $W$ is then uniform on every compact set of $\R_+$, and thus uniform on $\R_+$ since the functions are decresaing and bounded from below.
\end{demo}

The Laplace exponent $\psi$ of $Z$ is given for all $\lambda\geq0$ by : 
$$\psi(\lambda):=c\lambda+\frac12 b^2\lambda^2-\int (1-e^{-\lambda u}-\lambda h(u))\L(\d u),$$ 
where $h$ is a truncation function on $\R$ (see Section \ref{sec_SPLP}). Recall that $c$ depends on the choice of $h$. Then we have 

\begin{proposition}\label{prop_Cb_Ku2}
 Let $(g_n)_{n\geq0}$ and $g$ be continuous bounded mappings from $\R_+$ to $\R$, where $g$ satisfies $g(u)/u^2 \to K$ as $u\to 0+$ for some constant $K$. Assume that the mappings $\tilde g_n:u\mapsto \frac{g_n(u)}{1\wedge u^2}$ converge uniformly to  $\tilde g:u\mapsto \frac{g(u)}{1\wedge u^2}$ on $\R_+^*$. Then as $n\to\infty$,
\upshape $$\tL_n(g_n) \underset{n\to\infty}{\to} \L(g) + Kb^2  .$$
\end{proposition}

We first prove the following two lemmas. Define $\mathcal M_{\textsc L}(\R)$ the set of $\sigma$-finite measures $\nu$ on $(\R,\mathcal B(\R))$ satisfying the condition $\int (1\wedge |u|^2) \nu(\d u)<\infty$. 
\begin{lemme}\label{lemme_Cb_Ku2}
 Let $(h_n)_{n\geq0}$ and $h$ be continuous bounded mappings from $\R$ to $\R$, where $h$ satisfies $h(u)/u^2 \to K$ as $u\to 0$ for some constant $K$. Consider $(\nu_n)_{n\geq0}$ and $\nu$ in $\mathcal M_{\textsc L}(\R)$ and assume that :
\begin{enumerate}[\upshape(i)]
 \item There exists $a\in\R$ such that for all continuous bounded function $f$ satisfying $f(u)/u^2\to K$ as $u\to0$, 
\upshape$$\nu_n(f)\underset{n\to\infty}{\to} \nu(f)+K a.$$
\itshape
 \item The mappings $\tilde h_n:u\mapsto \frac{h_n(u)}{1\wedge u^2}$ converge uniformly to  $\tilde h:u\mapsto \frac{h(u)}{1\wedge u^2}$ on $\R^*$.
\end{enumerate}
Then \upshape $$\nu_n(h_n) \underset{n\to\infty}{\to} \nu(h) + K a.$$
\end{lemme}

\begin{demo}
First note that since $\nu_n,\nu\in\mathcal M_{\textsc L}(\R)$, all the integrals considered in the statement of the theorem are finite. Write :
$$  \left|\int h_n\d\nu_n-\int h\d\nu -Ka\right|\leq \left|\int (h_n-h)\d\nu_n\right| + \left|\int h\d\nu_n-\int h\d\nu-Ka\right|.$$
The mapping $h$ is continuous and bounded on $\R$, and satisfies $h(u)\sim K u^2$ when $|u|\to0$ ; then (i) implies the convergence to $0$ of the term $|\int h\d\nu_n-\int h\d\nu-K a|$. \\
Let $\eps$ be a positive real number. First observe that (ii) implies that $\tilde h_n$ and $\tilde h$ can be extended to continuous functions on $\R$ (which we will also denote by $\tilde h_n$ and $\tilde h$), and we have $\tilde h_n(0)\to\tilde h(0)=K$. Then (ii) implies for $n$ large enough and any $u\in\R$ : $$|\tilde h_n-\tilde h|(u)\leq \eps,$$ and then we have $|\int (h_n-h)\d\nu_n|\leq \eps \int (1\wedge u^2) \nu_n(\d u)$. Now according to (i), the sequence $(\int (1\wedge u^2) \nu_n(\d u))_{n}$ converges and is consequently bounded. This proves that $|\int (h_n-h)\d\nu_n|$ tends to $0$ and ends the proof.\end{demo}

\begin{lemme} \label{JS_Lambda_Ku2}
 Let $g$ be a continuous bounded function on $\R_+$ such that for some $K\in\R$, $g(u)/u^2 \to K$ as $u\to0$. Then
\upshape$$\int g\d\tL_n \to Kb^2+\int g\d\L \ \ \ \ when\ n\to\infty.$$
\end{lemme}
\begin{demo}
Considering Assumption \HypLevy, first notice that a straightforward application of Proposition \ref{prop_thJS} yields
\begin{enumerate}[\upshape(a)]
 \item For all truncation function $h$ on $\R_+$, $\int h^2\d\tL_n \to b^2+\int h^2\d\L$ as $n\to\infty$.
 \item For any continuous bounded function $g$ such that $g(u)=o(u^2)$ as $u\to0$, $\int g\d\tL_n \to \int g\d\L$.
\end{enumerate}
Then, let $h$ be a truncation function on $\R_+$. Writing $g=Kh^2+(g-Kh^2)$, we get :
\begin{align*}
\left|\int g\d\tL_n - (Kb^2+\int g\d\L)\right|\ \leq\ &\left|\int Kh^2\d\tL_n -(Kb^2+\int Kh^2\d\L)\right|\\&+\left|\int (g-Kh^2)\d\tL_n - \int(g-Kh^2)\d\L\right|. 
\end{align*}
Now since $h$ is a truncation function, thanks to (a) we know that the first term of the right-hand side vanishes as $n\to\infty$. As for the second term, the function $g-Kh^2$ is bounded and satisfies $\underset{u\to0}{\lim} \frac{g(u)-Kh^2(u)}{u^2}=0$ so that we can apply (b), and $|\int (g-Kh^2)\d\tL_n - \int(g-Kh^2)\d\L|\to0$ as $n\to\infty$.
\end{demo}

Finally, Proposition \ref{prop_Cb_Ku2} arises as a direct consequence of Lemmas \ref{lemme_Cb_Ku2} and \ref{JS_Lambda_Ku2}.

\paragraph*{Convergence of the classical ladder process}\quad

For all $n\geq1$ let $\kappa_n$ be the Laplace exponent of the bivariate ladder process $(\Linv_n, \Han)$, and denote by $\kappa$ the Laplace exponent of $(\Linv,\Ha)$. Note that the condition of normalization (\ref{formula_norm_local_time}) imposed to $L$ implies $\kappa(1,0)=\phi(1)^{-1}$. 

\begin{proposition} \label{prop_cv_H1}
The sequence $(\Linv_n, \Han)_{n\geq1}$ converges weakly in distribution to $(\Linv,\Ha)$.
\end{proposition}

\begin{lemme}\label{lemme_kappa1}
For all $n\geq1$, $\kappa_n(1,0)=\tphi_n(1)^{-1}$.
\end{lemme}
\begin{demo}
Let $T_n$ be the first jump time of the process $\bar \tZ_n$. 
The subordinator $\Linv_n$ is a compound Poisson process with rate $\alpha_n$ and jump size distribution $\mathcal L(T_n)$ (where $\mathcal L(T_n)$ denotes the law of $T_n$). Therefore we have
\begin{align*}
 e^{-\kappa_n(1,0)} &=\E(e^{-\Linv_n(1)})\\
		&=\sum_{k\geq0} \frac{(\alpha_n)^ke^{-\alpha_n}}{k!}\E(e^{-T_n})^k \\
		&=e^{-\alpha_n(1-\E(e^{-T_n}))}.
\end{align*}
Now the variable $T_n$ is a.s. finite and from Theorem 8.1 and Lemma 8.6 in \cite{K}, we get $$\E(e^{-T_n})=1-\frac{n}{d_n\tphi_n(1)}.$$ 
Since $\alpha_n=\frac{d_n}{n}$, we get $\alpha_n(1-\E(e^{-T_n}))=\tphi_n(1)^{-1}$ and in consequence $\kappa_n(1,0)=\tphi_n(1)^{-1}$.
\end{demo}

\begin{demopr}{\ref{prop_cv_H1}}
According to \cite[Th. VII.3.4]{JS}, proving the convergence of the Laplace exponents of $(\Linv_n,\Han)$ is sufficient. Fix $(\a,\b)\in\R_+\times\R_+$. From Corollary VI.10 in \cite{B}, and since $\tZ_n$ (resp. $Z$) is not a compound Poisson process (implying its marginal distributions do not have an atom at zero), we know that 
\[\kappa_n(\a,\b)=\kappa_n(1,0)\exp\Big\{\int_0^\infty \d t \int_{(0,\infty)} (e^{-t}-e^{-\a t-\b x})\frac{1}{t}\P(\tZ_n(t)\in\d x)\Big\} \]
and
\[\kappa(\a,\b)=\kappa(1,0)\exp\Big\{\int_0^\infty \d t \int_{(0,\infty)} (e^{-t}-e^{-\a t-\b x})\frac{1}{t}\P(Z(t)\in\d x)\Big\} .\]

First assume that $\beta=0$ and $\alpha>1$. From Assumption \HypLevy, for all $t>0$ a.s. the measures $\P(\tZ_n(t)\in\d x)\mathds1_{x>0}$ converge weakly to $\P(Z(t)\in\d x)\mathds1_{x>0}$. Besides, Lemma \ref{lemme_kappa1} ensures the convergence of $\kappa_n(1,0)=\tphi_n(1)^{-1}$ towards $\phi(1)^{-1}=\kappa(1,0)$ as $n\to\infty$. Then using Fatou's Lemma we obtain 
$$\liminf \kappa_n(\a,0)\geq \kappa(\alpha,0).$$

But from (3), p. 166 in \cite{B}, $\kappa_n(\a,0)\hat\kappa_n(\a,0)=\a$, so that 
$$\liminf \frac{1}{\hat\kappa_n(\a,0)}\geq \frac{1}{\hat\kappa(\a,0)},\ \text{and then }\ \limsup \hat\kappa_n(\a,0)\leq\hat\kappa(\a,0),$$
where $\hat\kappa_n$ and $\hat\kappa$ refer respectively to $\hat\tZ_n=-\tZ_n$ and $\hat Z=-Z$. Then replacing $\tZ_n$ by $\hat\tZ_n$ in the above arguments entails $\kappa_n(\a,0)\to\kappa(\a,0)$ as $n\to\infty$. The same arguments hold for $\alpha\in(0,1)$ by exchanging $\limsup$ and $\liminf$.

Now using the notation of Chapter VI in \cite{B}, let $\tau$ be an indepedent exponential variable with parameter $q>0$, and define $\Gtaun:=\sup\{t<\tau, \bar\tZ_n(t)=\tZ_n(t)\}$ (resp. $\Gtau:=\sup\{t<\tau, \bar Z(t)=Z(t)\}$) the last zero of the reflected process $\bar\tZ_n-\tZ_n$ (resp. $\bar Z-Z$) before $\tau$. \\

\textbf{Claim :}\qquad  $(\Gtaun,\bar\tZ_n(\tau))$ converges in law towards $(\Gtau,\bar Z(\tau))$.

\begin{description}
 \item [\qquad] From the weak convergence of $\tZ_n$ towards $Z$  and using the Skorokhod representation theorem, we can assume that $\tZ_n\xrightarrow{a.s.} Z$, and it is sufficient to prove the a.s. convergence of $(\Gtaun,\bar\tZ_n(\tau))$ towards $(\Gtau,\bar Z(\tau))$. \\

First, the a.s. convergence of $\bar\tZ_n(\tau)$ towards $\bar Z(\tau)$ is straightforward from Proposition VI.2.11 in \cite{JS}. Let us now prove that $\Gtaun\xrightarrow{a.s.}\Gtau$.

Using time reversal and considering the infimum process of the reversed reflected process, a direct adaptation of the proofs of Propositions VI.2.4 and VI.2.11 in \cite{JS} allows us to obtain the inequality $\liminf \Gtaun\geq \Gtau$ a.s.

Let us now prove that $\limsup \Gtaun\leq \Gtau$ a.s. Consider a realization $(\tilde z_n,z,g_T^{(n)},\gtau,T)$ of $(\tZ_n,Z,\Gtaun,\Gtau,\tau)$. Assume there exists $t\in(0,T)$, and a subsequence $(k_n)_{n\in\N}$ satisfying for all $n$, $\gtaukn<t<\gtau$.  

For all $n\in\N$, by definition of $\gtaukn$, on $[\gtaukn,T)$ the current supremum $\bar{\tilde z}_{k_n}$ of $\tilde z_{k_n}$ is equal to a constant $s^{(k_n)}$. Define $\bar z$ the current supremum of $z$ and  $g'_T:=\sup\{u<\gtau,\ \bar z(u)-z(u)=0\}$ the penultimate zero of $\bar z-z$ before $T$. On $[g'_T,\gtau)$ (resp. on $[\gtau,T)$), $\bar z$ is equal to a constant $s'$ (resp. $s$). Applying Proposition VI.2.11 in \cite{JS} at times $t$ and $g_T$, we obtain the convergence of $s^{(k_n)}$ towards $s$ and $s'$, which entails $s=s'$. Hence we get the existence of two times $g'_T<g_T$ such that $\bar z(g'_T)=z(g'_T)=\bar z(\gtau)=z(\gtau)$. 
Finally, Proposition VI.4 in \cite{B} shows that such realizations form a negligible set, so that $\limsup \Gtaun\leq \Gtau$ a.s. We conclude that $\Gtaun\xrightarrow{a.s.}\Gtau$.
\end{description}

The convergence in law of $(\Gtaun,\bar\tZ_n(\tau))$ towards $(\Gtau,\bar Z(\tau))$ entails, from $(1)$ p.163 in \cite{B}, that $\kappa_n(q,0)/\kappa_n(\a+q,\b)\to\kappa(q,0)/\kappa(\a+q,\b)$ as $n\to\infty$. We conclude from the convergence established above.
% On the other hand, thanks to the lemma above, we know that for all $n\geq1$, $\kappa_n(1,0)=\tphi_n(1)^{-1}$, which converges as $n\to\infty$ towards $\phi(1)^{-1}=\kappa(1,0)$, and we have the announced convergence.
\end{demopr}
\\

\paragraph{Proof of Theorem \ref{th_cv_H}} \quad

The proof of the theorem will consist in applying Proposition \ref{prop_thJS} to the sequence of bivariate Lévy processes $(\Han,\Hcn)$. To this aim we first establish the following property.

\begin{proposition} \label{prop_mu1}
The measure $\mu(\,\cdot\,,\un)$ is finite, and for any continuous bounded function $g$ on $\R_+$ which is differentiable at 0, we have as $n\to\infty$ :\\
\begin{enumerate}[\upshape(i)]
 \item \label{mu1}
Under Assumption \HypMutConst, \qquad \upshape$\int g(u)\mu_n(\d u,\un)\to  \theta g(0).$\itshape\\
Under Assumption \HypMut, \qquad \upshape$\int g(u)\mu_n(\d u,\un)\to \int g(u)\mu(\d u,\un)+\rho g(0).$\\
\itshape where $\rho=\kappa b^2$ has been defined in Theorem \ref{th_cv_H}.\\

 \item \label{mu0}
Now if $g(0)=0$, \\
Under Assumption \HypMutConst, \qquad \upshape $\int g(u)\mu_n(\d u,\zero)\to \int g(u)\mu^+(\d u)+g'(0)\frac {b^2}2.$ \itshape\\
Under Assumption \HypMut, \qquad \upshape$\int g(u)\mu_n(\d u,\zero)\to \int g(u)\mu(\d u,\zero)+g'(0)\frac {b^2}2.$\\

\end{enumerate}
\itshape Furthermore, in both cases, for all $\delta>0$, the results are still valid if we replace $g$ by $g\,\mathds1_{[0,\delta]}$ or by $g\,\mathds1_{(\delta,\infty)}$.
\end{proposition}

First of all, to prove this proposition we need the two lemmas below. The first one is deduced from the convergence in law of $(\Ha_n)$. The second one is specific to the case \HypMutConst.

\begin{lemme} \label{lemme_cv_H1}
Let $g$ be a continuous bounded function from $\R_+$ to $\R$ such that $g(u)=o(u^2)$ as $u\to0$. We have
\upshape $$\int_{(0,\infty)} \left( \int_0^z g(z-y)e^{-\teta_n y} \d y \right)\tL_n(\d z) \underset{n\to\infty}\to  \int_{(0,\infty)} \left( \int_0^z g(z-y)e^{-\eta y} \d y \right) \L(\d z).$$
\end{lemme}

\begin{demo}
From Proposition \ref{prop_thJS} and the convergence $\Han\Rightarrow\Ha$ established in Proposition \ref{prop_cv_H1}, we get that $\mu^+_n(g) \to\mu^+(g)$ as $n\to\infty$, where $\mu^+$ denotes the Lévy measure of $\Ha$. Now we deduce from the expression of $\mu^+_n$ given by (\ref{mu+_n}) that 
$$ \mu^+_n(g) = \int_{(0,\infty)} \tL_n(\d z) \int_0^z e^{-\teta_n y} g(z-y)\d y.$$
A similar calculation for the limiting process gives $\mu^+(g) = \int_{(0,\infty)} \L(\d z) \int_0^z e^{-\eta y} g(z-y)\d y$, and the result follows.
\end{demo}

\begin{lemme} \label{lemme_theta_n}
 As $n\to\infty$, we have \upshape
$$\int_{(0,\infty)} \tL_n(\d u)\left(\int_0^{1\wedge u} e^{\teta_n (r-u)}\d r\right) \sim \frac{d_n}{n}.$$
\end{lemme}

\begin{demo}
 For all $a>0$, we have by definition of $\tilde\phi_n$ and thanks to formula \eqref{formula_LK_vf}:
$$\frac{d_n}{n}\tphi_n(a)-\int_{(0,\infty)}(1-e^{-\tphi_n(a)u})\tL_n(\d u) = a.$$
Then we have
$$1-\frac{n}{d_n}\int_{(0,\infty)}\frac{1-e^{-\tphi_n(a) u}}{\tphi_n(a)} \tL_n(\d u)=\frac{n}{d_n}\frac{a}{\tphi_n(a)}$$
which leads to
\begin{multline*}
\frac{n}{d_n}\int_{(0,\infty)}\tL_n(\d u)\left(\int_0^{1\wedge u} e^{\teta_n (r-u)}\d r \right) \\
=1-\frac{n}{d_n}\frac{a}{\tilde\phi_n(a)}-\frac{n}{d_n}\int_{(0,\infty)}\left(\frac{1-e^{-\tilde\phi_n(a) u}}{\tilde\phi_n(a)}-\int_0^{1\wedge u} e^{\teta_n (r-u)}\d r\right) \tL_n(\d u) .
\end{multline*}
Now it is easy to check that we can apply Proposition \ref{prop_Cb_Ku2} (further applications of this proposition are detailed in the proof of Proposition \ref{prop_mu1}) to get the convergence of 
$$\int_{(0,\infty)}\left(\frac{1-e^{-\tphi_n(a) u}}{\tphi_n(a)}-\int_0^{1\wedge u} e^{\teta_n (r-u)}\d r\right) \tL_n(\d u) $$ towards a finite quantity. Furthermore, we know that $\tilde\phi_n(a)\to\phi(a)$ and that $\frac{n}{d_n}$ vanishes as $n\to\infty$, which leads to the announced result. 
\end{demo}

\par\bigskip

\begin{demopr}{\ref{prop_mu1}}
We begin with the proof of point \eqref{mu1}. Let $g$ be a continuous bounded function on $\R_+$, differentiable at $0$. We have : 
\begin{align*}
 \int g(u) \mu_n(\d u,\un)
&=\int_{(0,\infty)} \int_0^\infty \d y\; e^{-\teta_n y}\; \tL_n(y+\d u)\; f_n(n(u+y))\; g(u) \\
&=\int_{(0,\infty)} \tL_n(\d u)\; f_n(nu)\;\int_0^u \;\d y\; e^{-\teta_n y}\; g(u-y) \\
&=\int_{(0,\infty)} \tL_n(\d u)\; f_n(nu) \;\int_0^u \;\d z\; e^{\teta_n (z-u)}\; g(z), 
\end{align*}
and a similar calculation is available for $\mu$.\\
Let us first treat the case of Assumption \HypMut. The calculation above entails
\begin{multline}\label{inequation}
 \Big|\int g(u) \mu_n(\d u,\un)-\int g(u) \mu(\d u,\un)-\rho g(0)\Big| \\
			 \leq \quad \Big|\int \tL_n(\d u)f_n(nu)\bigg(\int_0^{1\wedge u} \d z\, e^{\teta_n (z-u)} g(z)\bigg)-\int \L(\d u)f(u)\bigg(\int_0^{1\wedge u} \d z\, e^{\eta (z-u)} g(z)\bigg)  -\rho g(0)\Big|\\
			  +\Big|\int \tL_n(\d u) f_n(nu)\bigg(\int_{1\wedge u}^u \d z\, e^{\teta_n (z-u)} g(z)\bigg) -\int \L(\d u)f(u)\bigg(\int_{1\wedge u}^u \d z\, e^{\eta (z-u)} g(z)\bigg)\Big|.
\end{multline}
First note that the integral $\int_{1\wedge u}^u \d z\, e^{\teta_n (z-u)} g(z)$ can be rewritten as $\int_0^u \d y\, e^{-\teta_n y} g(u-y) \mathds1_{u-y\geq 1}$. Since the function $z\mapsto g(z) \mathds1_{z\geq 1}$ is bounded and vanishes on $[0,1]$, a simple approximation argument allows us to obtain from Lemma \ref{lemme_cv_H1} the convergence of $\int \tL_n(\d u) \big(\int_{1\wedge u}^u \d z\, e^{\teta_n (z-u)} g(z)\big)$ towards $\int \L(\d u)\big(\int_{1\wedge u}^u \d z\, e^{\eta (z-u)} g(z)\big)$. Then, the convergence to $0$ of the second term in the right-hand side is obtained using the fact that $|f_n|\leq1$ for all $n$, and the uniform convergence on $\R_+$ of $f_n(n\cdot)$ towards $f$. \\

Next we focus on the first term. We set $h_n(u):=f_n(nu)\int_0^{1\wedge u} \d z\, e^{\teta_n (z-u)} g(z)$ and $\ h(u):=f(u)\int_0^{1\wedge u} \d z\, e^{\eta (z-u)} g(z)$. The aim of the next paragraph is to check that the functions $h_n$ and $h$ satisfy the hypotheses of Proposition \ref{prop_Cb_Ku2}, which will entail the convergence to $0$ of the first term in the right-hand side of (\ref{inequation}).

\begin{itemize}
\item The functions $|h_n|$ and $|h|$ can be upper bounded by $\int_0^1 |g(z)|\d z$, which is a finite quantity. Moreover the continuity of $g$, $f_n$ and $f$ ensures that of $h_n$ and $h$.
\item We have for $u\leq1$
$$\frac{f(u)}{u}\times\underset{x\in[0,u]}{\min} \{g(x)\} \frac1u\int_0^u e^{-\eta y}\d y \leq \frac{h(u)}{u^2}\leq \frac{f(u)}{u}\times\underset{x\in[0,u]}{\max} \{g(x)\} \frac1u\int_0^u e^{-\eta y}\d y,$$
Now $\frac1u \int_0^u e^{-\eta y}\d y \to 1$ as $u\to0$, and $\underset{u\to0}\lim \underset{x\in[0,u]}{\min} \{g(x)\}=\underset{u\to0}\lim \underset{x\in[0,u]}{\max} \{g(x)\}=g(0)$ (recall that $g$ is continuous). Then thanks to Assumption \HypMut.\eqref{HypMut_0} we can conclude that $\underset{u\to0}{\lim} \frac{h(u)}{u^2}=\kappa g(0)$. Besides, this conclusion ensures that $\mu(\,\cdot\,,\un)$ is a finite measure.
 
\item Finally, the mappings $u\mapsto \frac{h_n(u)}{1\wedge u^2}$ converge to $u\mapsto \frac{h(u)}{1\wedge u^2}$ uniformly on $\R^*$. 
Indeed, fix $\eps>0$. For all $u\in(0,1)$,
\begin{align*}
 \bigg|\frac{h_n(u)-h(u)}{u^2}\bigg|
&\leq\frac1{u^2}\underset{[0,1]}{\max}|g|\int_0^u \Big(|f_n(nu)-f(u)|e^{-\teta_n y}+f(u)|e^{-\teta_n y}-e^{-\eta y}| \Big)\;\d y \\
&\leq\underset{[0,1]}{\max}|g|\,\bigg(\frac{|f_n(nu)-f(u)|}{u}\frac1u\int_0^u \d y+|\teta_n-\eta|\frac{1}{u^2}\int_0^u y\;\d y\bigg)\\
&\leq\underset{[0,1]}{\max}|g|\,\bigg(\frac{|f_n(nu)-f(u)|}{u}+\frac12 |\teta_n-\eta|\bigg).
\end{align*}
Then thanks to Assumption \HypMut.\eqref{HypMut_cv}, and since $\teta_n\to\eta$, for $n$ large enough $\big|\frac{h_n(u)-h(u)}{u^2}\big|\leq\eps$ for all $u\in(0,1)$. \\
On the other hand, for all $u\geq1$, 
$$
 |h_n(u)-h(u)| 
\leq \underset{[0,1]}{\max}|g|\,(|f_n(nu)-f(u)|+\frac12 |\teta_n-\eta|)
,$$
which can be upper bounded by $\eps$ for all $u\geq 1$ and $n$ large enough, again thanks to Assumption \HypMut.\eqref{HypMut_cv} and to the convergence of $\teta_n$ towards $\eta$. \\
It follows then that for all $u\in\R^*$ and $n$ large enough, $\big|\frac{h_n(u)}{1\wedge u^2}\big|\leq \eps$.\\
\end{itemize}
All the conditions of Proposition \ref{prop_Cb_Ku2} are then fulfilled, and we get the claimed convergence.\\

We now consider Assumption \HypMutConst. Exactly as before, we have
\begin{align*}
 \Big|\int g(u) \mu_n(\d u,\un)-\theta g(0)\Big| &\leq \quad \Big|\theta_n \int \tL_n(\d u) \bigg(\int_0^{1\wedge u} \d z\, e^{\teta_n (z-u)} g(z)\bigg)-\theta g(0)\Big|\\
			  &\ \ +\Big|\theta_n \int \tL_n(\d u)\bigg(\int_{1\wedge u}^u \d z\, e^{\teta_n (z-u)} g(z)\bigg) \Big|,
\end{align*}
and as in case \HypMut, Lemma \ref{lemme_cv_H1} entails the convergence to $0$ of the second term in the right-hand side. \\

As for the first term, we have
\begin{align*}
 \Big|\theta_n  \int \tL_n(\d u)\bigg(\int_0^{1\wedge u} \d z\, e^{\teta_n (z-u)}& g(z)\bigg)-\theta g(0)\Big|\\
 &\leq \theta_n \int \tL_n(\d u) \bigg(\int_0^{1\wedge u} \d z\, e^{\teta_n (z-u)} |g(z)-g(0)|\bigg)\\
 &+\Big|\theta_n \int \tL_n(\d u) \bigg(\int_0^{1\wedge u} \d z\, e^{\teta_n (z-u)} g(0)\bigg)-\theta g(0)\Big|,
\end{align*}
Lemma \ref{lemme_theta_n} ensures the convergence to $0$ of the second term in the right-hand side. Now the functions $u\mapsto \int_0^{1\wedge u} \d z\, e^{\teta_n (z-u)} |g(z)-g(0)|$ are continuous, bounded by $\sup g$, and converge to $u\mapsto \int_0^{1\wedge u} \d z\, e^{\eta (z-u)} |g(z)-g(0)|$, which is equivalent to $g'(0)u^2/2$ as $u\to0$. As a consequence of Proposition \ref{prop_Cb_Ku2}, we then have the convergence of $\int \tL_n(\d u) \big(\int_0^{1\wedge u} \d z\, e^{\teta_n (z-u)} |g(z)-g(0)|\big) $ to
$g'(0)\frac{b^2}2+\int \tL(\d u) \big(\int_0^{1\wedge u} \d z\, e^{\eta (z-u)} |g(z)-g(0)|\big)$, which is a finite quantity, 
and thus the fact that $\theta_n\to0$ ends the proof of the second assertion in \eqref{mu1}.\\
 
The proof of point \eqref{mu0} is very similar :
Under Assumption \HypMutConst\ or \HypMut, we have 
$$\int g(u) \mu_n(\d u,\zero)=\int h_n(u) \;\tL_n(\d u)$$ 
with $$h_n(u):=(1-f_n(nu))\;\int_0^u \;\d z\; e^{\teta_n (z-u)}\; g(z).$$ 
The same arguments as in the proof above work, except for the limit at $0$ of $h(u)/u^2$ : in this case, the fact that $g(u)/u\to g'(0)$ as $u\to0$ implies $\frac1{u^2}\int_0^u e^{-\eta y}g(u-y)\d y\to \frac {g'(0)}2$, and then since $1-f(u)\to1$, we get $\frac{h(u)}{u^2}\to\frac {g'(0)}2$ as $u\to0$. Finally, in the case of Assumption \HypMutConst, $f\equiv0$ implies $\mu(\d u,\zero)=\mu^+(\d u)$, which allows us to conclude.\\

To get the last conclusion of the proposition, first notice that $\mu$ has no atom : Suppose $\mu$ has an atom $d>0$, then $\mu(\{d\})=\int_0^\infty e^{-\eta x}\L(\{x+d\}) \d x>0$, which leads to the existence of a subset $U\subset[d,+\infty)$ such that Leb$(U)\neq0$ and $\L(\{y\})>0$ for all $y\in U$. This implies $\L(U)=+\infty$, which is impossible since $\L(U)\leq\L([d,+\infty))<\infty$.\\ 
The results follow then by approximation : for all $\eps>0$, let $I_\eps^+$ and $I_\eps^-$ be two continuous piecewise linear functions satisfying :
\[
I_\eps^+(x) = \left \{ \begin{array}{lrl}
			0 & \text{ if} & x \leq \delta \\
			1 & \text{ if} & x \geq \delta+\eps \\
\end{array}
\right.
\qquad
I_\eps^-(x) = \left \{ \begin{array}{lrl}
			0 & \text{ if} & x \leq \delta-\eps \\
			1 & \text{ if} & x \geq \delta \\
\end{array}\right. .
\]
We have $I_\eps^-\leq\mathds1_{[0,\delta]}\leq I_\eps^+$. This gives, for all $\eps>0$,
$$\int gI_\eps^-\d\mu+\rho g(0) \leq \underset{n\to\infty}\liminf \int_{[0,\delta]} g\d\mu_n \leq \underset{n\to\infty}\limsup \int_{[0,\delta]} g\d\mu_n \leq \int gI_\eps^+\d\mu+\rho g(0).$$
Now when $\eps\to0$, $\int gI_\eps^-\d\mu\to\int_{[0,\delta]} g\d\mu$ and $\int gI_\eps^-\d\mu\to\int_{[0,\delta)} g\d\mu$. Since $\mu$ has no atom, these two integrals are equal and we get 
$$\int_{[0,\delta]} g(u)\mu_n(\d u,\un)\to \int_{[0,\delta]} g(u)\mu(\d u,\un)+\rho g(0) .$$ 
The other announced results can be obtained by a similar reasoning.
\end{demopr}

\par\medskip
\begin{demoth}{\ref{th_cv_H}}
We first prove the second part of the theorem, i.e. we assume \HypMut. Moreover we assume first that $Z$ does not drift to $-\infty$. Proposition \ref{prop_mu1} allows us to establish the three claims below, which correspond respectively to points (\ref{thJS_mesure_Levy}), (\ref{thJS_drift}) and (\ref{thJS_coeff_brownien}) of Proposition \ref{prop_thJS}.\\

\begin{claim}
For all continuous bounded function $g$ on $\R_+^2$ such that $g$ is zero in a neighbourhood of $(0,0)$, 
$$ \mu_n(g)\to  (\mu+\rho\delta_{(0,1)})(g).$$
\end{claim}
We have :
\begin{itemize}
 \item First, since $u\mapsto g(u,0)$ is zero in a neighbourhood of 0, $$\int g(u,0) \mu_n(\d u,\zero) \to \int g(u,0) \mu(\d u,\zero)$$ as $n\to\infty$ thanks to Proposition \ref{prop_mu1} \eqref{mu0}.
 \item Second $\int g(u,1) \mu_n(\d u,\un) \to \int g(u,1) \mu(\d u,\un) + \rho g(0,1)$ according to Proposition \ref{prop_mu1} \eqref{mu1}.
\end{itemize}
and the result follows.\\

\begin{claim}
 For all $(\alpha,\beta)\in\R_+^2$, 
\upshape $$\int (\alpha,\beta)\tr h(u,q) \mu_n(\d u,\d q) \to \frac{b^2}2 \alpha + \rho\beta +\int (\alpha,\beta)\tr h(u,q) \mu(\d u,\d q),$$
where $h$ is the truncation function defined earlier.
\end{claim}

We have $$
 \int (\alpha,\beta)\tr h(u,q) \mu_n(\d u,\d q) 
=\int_{[0,\delta]} (\alpha u+\beta q) \mu_n(\d u,\d q) +\int_{(\delta,\infty)} (\alpha\delta+\beta q) \mu_n(\d u,\d q),$$
and then :
\begin{itemize}
 \item $u\mapsto\alpha u+\beta$ is a continuous bounded function on $[0,\delta]$, then thanks to Proposition \ref{prop_mu1}, 
$$\int_{[0,\delta]} (\alpha u+\beta) \mu_n(\d u,\un)\to \rho\beta+\int_{[0,\delta]} (\alpha u+\beta) \mu(\d u,\un).$$
 \item In the same way, thanks to Proposition \ref{prop_mu1} \eqref{mu0}, $$\int_{[0,\delta]} \alpha u \mu_n(\d u,\zero)\to \frac{b^2}2 \alpha+\int_{[0,\delta]} \alpha u\mu(\d u,\zero).$$
 \item And finally, thanks to Proposition \ref{prop_mu1} (points \eqref{mu1} and \eqref{mu0}), 
$$\int_{(\delta,\infty)} (\alpha\delta+\beta q) \mu_n(\d u,\d q)\to\int_{(\delta,\infty)} (\alpha\delta+\beta q) \mu(\d u,\d q)\ \ \text{ when } n\to\infty.$$
\end{itemize}
As a consequence, 
\begin{align*}
&\int (\alpha,\beta)\tr h(u,q) \mu_n(\d u,\d q)\\
&\to \rho\beta+\int_{[0,\delta]} (\alpha u+\beta) \mu(\d u,\un)+\frac{b^2}2\alpha+\int_{[0,\delta]} \alpha u\mu(\d u,\zero)+\int_{(\delta,\infty)} (\alpha\delta+\beta q) \mu(\d u,\d q) \\
&= \rho\beta + \frac{b^2}2\alpha +\int(\alpha,\beta)\tr h(u,q) \mu(\d u,\d q),
\end{align*}
which proves our assertion.\\

\begin{claim}
Denote by $h_1$ (resp. $h_2$) the first (resp. second) coordinate of $h$. For all $i,j\in\{1,2\}$, \upshape
$$\int h_i(u,q)h_j(u,q)\mu_n(\d u,\d q) \underset{n\to\infty}{\to} \int h_i(u,q)h_j(u,q)(\mu(\d u,\d q)+\rho \delta_0(\d u)\delta_1(\d q))$$
\itshape as $n\to\infty$.
\end{claim}

Note that $\int h_i(u,q)h_j(u,q)\delta_0(\d u)\delta_1(\d q)=h_i(0,1)h_j(0,1)$.
 \begin{itemize}
  \item The continuous bounded function $h_1^2$ satisfies $h_1(u,q)^2/u\to0$ as $u\to0$, for $q\in\{0,1\}$. Then thanks to Proposition \ref{prop_mu1} (points \eqref{mu1} and \eqref{mu0}) we have 
$$\int h_1(u,q)^2\mu_n(\d u,\d q) \underset{n\to\infty}{\to} \int h_1(u,q)^2 \mu(\d u,\d q),$$
and since $h_1(0,1)=0$, we get the announced result for $(i,j)=(1,1)$.
 \item The continuous bounded function $u\mapsto h_1(u,1)h_2(u,1)$ satisfies $h_1(0,1)h_2(0,1)=0$ as $u\to0$, so that according to Proposition \ref{prop_mu1} \eqref{mu1}, 
$$\int h_1(u,1)h_2(u,1)\mu_n(\d u,\un) \underset{n\to\infty}{\to} \int h_1(u,1)h_2(u,1) \mu(\d u,\un).$$
Moreover, $h_1(0,1)=0$ and $h_2(u,0)=0$ for all $u\geq0$, and then we can deduce the result for $(i,j)=(1,2)$.
 \item Finally, when $q=0$ or $q=1$, we have $h_2(u,q)^2\equiv q$ for all $u\in\R_+$. In consequence, 
$$\int h_2(u,1)^2\mu_n(\d u,\un) \underset{n\to\infty}{\to} \int h_2(u,1)^2 \mu(\d u,\un)+\rho h_2(0,1)^2,$$
and since $h_2(u,0)\equiv0$, we get the result for $(i,j)=(2,2)$.
 \end{itemize}
\par\bigskip
Finally the three claims establish the theorem under Assumption \HypMut\ through a straightforward application of Proposition \ref{prop_thJS}. The proof in the case of Assumption \HypMutConst\ is very similar, and since in this case $f\equiv0$, the limiting Lévy measure is 
$$\mu(\d u,\{0,1\})+\theta\delta(\d q)\delta_0(\d u)=\mu^+(\d u)\delta_0(\d q)+\theta\delta_1(\d q)\delta_0(\d u),$$ 
which gives the expected result.\\

Finally we prove the theorem in the case where $Z$ drifts to $-\infty$. Using the convention that an exponentially distributed variable with parameter $0$ is equal to $+\infty$ a.s., and setting $\kill_n:=0$ when $\tZ_n$ does not drift to $-\infty$, all that is needed now is to prove that $\kill_n\to\kill$ as $n\to\infty$. Now since $W(\infty)<+\infty$, from the uniform convergence on $\R_+$ of $\tilde W_n$ towards $W$ (Proposition \ref{prop_cv}), we have $\tilde W_n(\infty)\to W(\infty)$, which ends the proof.
\end{demoth}

\section{Joint convergence in distribution of $\tZ_n$ with its local time at the supremum and its marked ladder height process} \label{sec_jointcv}

In this section we assume that Assumption \HypLevy, and one of the two Assumptions \HypMutConst\ or \HypMut\ hold, and we establish the joint convergence in law of $(\tZ_n,L_n,\Han,\Hcn)$. To prove this result, we will need the convergence in distribution of $\Hbn$ established in Section \ref{sec_cvMLHP}, and the joint convergence in distribution of $\tZ_n$ with its local time at the supremum and its classical ladder height process. The latter convergence is proved in L. Chaumont and R.A. Doney \cite{CD}, in the case of Lévy processes for which $0$ is regular for the open half-line $(0,\infty)$. We adapt here their proofs to our case of spectrally positive Lévy processes with finite variation.

\begin{theoreme} \label{th_appendix}
 The following convergence in distribution holds in $\D(\R)^4$ as $n\to\infty$ :
$$(\tZ_n,L_n,\Han,\Hcn)\Rightarrow(Z,L,\Ha,\Hc).$$
\end{theoreme}

This theorem is a consequence of the following proposition :

\begin{proposition} \label{prop_appendix}
 We have the following joint convergence in distribution in $\D(\R)^4$ as $n\to\infty$ :
$$(\tZ_n,L_n,\Han,\Hbn)\Rightarrow(Z,L,\Ha,\Hb).$$
\end{proposition}

\begin{demoth}{\ref{th_appendix}}
Consider the process $(\Han+\Hbn)$, denote by $\pi^\pm$ its jump point process (with values in $\R_+^*\times\{\partial\}$), and define $A:=\{t\in\R_+,\ \pi^\pm(t)\in\R_+^*\}$. Then we define the random process $\pi^\m$ as follows : conditional on $(\Han+\Hbn)$, for any $t$ in the countable set $A$, $\pi^\m(t)$ follows a Bernoulli distribution with parameter $\B{f_n(\pi^\pm(t))}$, and for $t\notin A$, $\pi^\m(t)=\partial$. Then by definition the process $\Hcn$ is distributed as a Poisson process with jump point process $\pi^\m$. It follows that conditional on $(\Han,\Hbn)$, the process $\Hcn$ is independent of $\tZ_n$ and $L_n$. Then Theorem \ref{th_cv_H} along with Proposition \ref{prop_appendix} entail the joint convergence in distribution of $(\tZ_n,L_n,\Han,\Hbn,\Hcn)$ towards $(Z,L,\Ha,\Hb,\Hc)$, and Theorem \ref{th_appendix} follows.
\end{demoth}\\

We now want to prove Proposition \ref{prop_appendix}, for which our inspiration comes from L. Chaumont and R.A. Doney \cite{CD}. With this aim in view, we need to introduce some notions about random walks. We consider the random walk $S=(S(j))_{j\geq0}$ defined by $S(0)=0$ and $S(j)=\sum_{i=1}^j Y_i$ for $j\geq1$, where $(Y_i)_{i\geq1}$ is a sequence of i.i.d. $\R$-valued random variables. We endow our random walk $S$ with a sequence of i.i.d. exponential random variables $(a_i)_{i\geq1}$ (their common parameter can be chosen arbitrarily), independent of $S$. We write $(N_t)_{t\geq0}$ for the Poisson process associated with this sequence of variables. We denote by $\bar S(j)$ the maximum of the random walk at step $j$ : $\bar S(j):=\max\{S(i),\ 1\leq i\leq j\}$, and we define its local time at the maximum :
$$\k(j):=\#\{i\in\{1,\ldots,j\},\ S(i)>\bar S(i-1)\}.$$
We then introduce a continuous-state version of the local time of $S$ at its maximum by setting
$$K(j):=\sum_{i=1}^{\k(j)} a_i.$$

We denote by $\t$ the right inverse of $\k$ :
$$\t(0)=0,\ \t(j+1)=\min\{i>\t(j),\ S(i)>S(\t(j))\},$$ 
which implies $\k(\t(j))=j$ for all integer $j$. Then similarly for $K$, we define $T$ by
$$\forall s\geq0,\ T(s)=\inf\{h\geq0, K(h)>s\},$$
which satisfies $T=\t\circ N$. Finally, we define $\g$ and $G$ as follows :
$$\forall j\geq0,\ \g(j)=\bar S(\t(j)),\ \text{ and } \forall s\geq0,\ G(s)=\bar S(T(s)).$$
The pair of processes $(\t,\g)$ is called ladder process, $\t$ being the ladder time process, and $\g$ the ladder height process. The pair $(T,G)$ is then a continuous-time version of the classical ladder process $(\t,\g)$.\\

In the sequel, we will consider a sequence of random walks $(S_n)_{n\geq1}$ (whose distributions can depend on $n$). As before, and independently for all $n$, we endow the random walk $S_n$ with a sequence of i.i.d. exponential variables $(A_i^n)_{i\geq1}$, independent of $S_n$, with parameter $\alpha_n$ to be specified later, and we denote by $N^n$ the corresponding Poisson process. We will use an obvious notation with subscript $n$ for all the quantities involved by the random walk $S_n$.\\

Let $X$ be a spectrally positive Lévy process (which is not a subordinator) with finite variation. We define its local time $L_X$ as in Section \ref{sec_Exc} : 
\[L_X(t):=\sum_{i=0}^{\ll(t)} A_i,\]
where $\ll(t)$ represents the number of jumps of the supremum until time $t$, and $(A_i)_{i\geq0}$ is a sequence of i.i.d. random exponential variables with arbitrarily chosen parameter $\alpha$, independent from $X$. We denote by $(\Linv_X,H)$ its bivariate ladder process and by $\kappa$ the Laplace exponent of the latter.

We define the convergence in distribution (resp. a.s.) of the sequence $(S_n)$ towards $X$ to be equivalent to the convergence in distribution (resp. a.s.) of the sequence of continuous-time processes $(S_n[nt])_{t\geq0}$ towards $X$, in $\D(\R_+)$. We keep again the notation $S_n\Rightarrow X$ for the convergence in law of $S_n$ to $X$.\\

The following four statements are the respective analogues of Theorem 1, Theorem 2, Theorem 3 and Corollary 2 in \cite{CD}, in the case of Lévy processes for which $0$ is not regular for the open half-line $(0,\infty)$. Our proofs are widely inspired of that of Chaumont and Doney in this paper.\\ 
\begin{proposition} \label{CD_th1} %analogue théorème 1
Let $(S_n)$ be a sequence of random walks converging in distribution to the Lévy process $X$. We then have the following convergence in law :
\upshape $$\left(\frac1n T_n,G_n\right) \Rightarrow (\Linv_X,H),$$
\itshape where for all $n$, the parameter $\alpha_n$ of the Poisson process $N^n$ is given by 
\upshape $$\alpha_n:=\exp\{\sum_{k\geq1} \frac1k e^{-k/n}\P(S_n(k)>0)\}.$$
\end{proposition}

\begin{demo}
The key of the following calculation is Fristedt's formula, which can be found in \cite[th. 10]{DoneyFluct} : 
$$1-\E(e^{-\delta \t_n(1) -\beta \g_n(1)})=\exp\{-\sum_{k\geq1} \frac{e^{-\delta k}}{k}\E(e^{-\beta S_n(k)},\ S_n(k)>0)\}.$$
It allows us to calculate the Laplace transform of $(\frac1n T_n,G_n)$ for all $\delta,\beta>0$ :
\begin{align*}
\E( & e^{-\delta T_n(1) -\beta G_n(1)})=\E( e^{-\delta \t_n(N_1^n) -\beta \g_n(N_1^n)})\\
&=\sum_{j\geq0} \E(e^{-\delta \t_n(1) -\beta \g_n(1)})^j\P(N_1^n=j)\\
&=e^{-\alpha_n} \sum_{j\geq0} \left(1-\exp\left\{-\int_{1/n}^\infty\frac{n}{[nt]}e^{-\delta [nt]/n} \E(e^{-\beta S_n([nt])},\ S_n([nt])>0)\d t\right\}\right)^j\frac{(\alpha_n)^j}{j!} \\
&=\exp\left\{ -\alpha_n \exp\left( -\int_{1/n}^\infty\frac{n}{[nt]}e^{-\delta [nt]/n} \E(e^{-\beta S_n([nt])},\ S_n([nt])>0)\d t \right)\right\}
.\end{align*}
Now from the expression of $\alpha_n$ we have
\begin{equation}\label{formula_alpha_n}
\alpha_n=\exp\left( \int_{1/n}^\infty \frac{n}{[nt]}e^{-[nt]/n}\P(S_n([nt])>0) \d t \right),
\end{equation}
and the convergence of $S_n$ towards $X$ gives, with an argument of dominated convergence as in the proof of Proposition \ref{prop_cv_H1},
\begin{align*}
 \alpha_n \exp & \left( -\int_{1/n}^\infty\frac{n}{[nt]}e^{-\delta [nt]/n} \E(e^{-\beta S_n([nt])},\ S_n([nt])>0)\d t \right) \\
& \underset{n\to\infty}\to \exp \left( -\int_0^\infty \left(\frac{e^{-t}}t \P(X_t>0)-\frac{e^{-\delta t}}t \E(e^{-\beta X_t},\ X_t>0) \right)\d t\right) \\
&\ \ \ = \exp \left( -\int_0^\infty \frac1t \E(e^{-t}-e^{-\delta t-\beta X_t},\ X_t>0)\d t\right) \\
&\ \ \ =\kappa(\delta,\beta),
\end{align*}
according to Corollary VI.10 in \cite{B}. \\
Thus we get the convergence of the Laplace exponent of $(\frac1n T_n,G_n)$ towards that of $(\Linv_X,H)$, which ends the proof.
\end{demo}

\begin{corollaire}\label{coroll_alpha}
 The parameters $\alpha_n$ converge to $\alpha$ as $n\to\infty$.
\end{corollaire}
\begin{demo}
 We saw in the proof above (see formula \eqref{formula_alpha_n} and following computation) that as $n\to\infty$, 
$$\alpha_n\to\exp\Big\{\int_0^\infty \frac{e^{-t}}t \P(X_t>0) \d t\Big\}.$$ 
Now this quantity is equal to $\kappa(\infty,0):=\underset{\delta\to\infty}\lim \kappa(\delta,0)$, and we have
$$\exp(\kappa(\infty,0))=\underset{\delta\to\infty}\lim \E(e^{-\delta \Linv_X(1)})=\P(A_1>1)=e^{-\alpha}.$$
\end{demo}
 
\begin{proposition} \label{CD_th2} %analogue théorème 2
 Under the same statement as in Proposition \ref{CD_th1}, assuming furthermore that the convergence of $(S_n)$ towards $X$ holds almost surely, for all fixed $t\geq0$ we have the convergence in probability of $K_n([nt])$ towards $L_X(t)$.
\end{proposition}

\begin{demo}
 Fix $\eps>0$ and $t\geq0$. Recall from the definitions of $K_n$ and $L_X$ that for all $n\geq1$, $j\geq0$, $t\geq0$,
$$K_n(j)=\sum_{i=1}^{\k_n(j)} A_i^n\ \text{ and }\ L_X(t)=\sum_{i=1}^{\ll(t)} A_i.$$
Write
$$\P\left(\Big|\sum_{i=1}^{\k_n([nt])} A_i^n - \sum_{i=1}^{\ll(t)} A_i\Big|>\eps\right) \leq \P\left(|\k_n([nt])-\ll(t)|>0\right)+ \P\left(\sum_{i=1}^{\ll(t)} |A_i-A_i^n|>\eps\right).$$

Fix $\eta>0$. On the one hand, since $\k_n$ and $\ll$ are finite integers, the almost sure convergence of $(S_n)$ towards $X$ ensures that for all $t\geq0$, $\k_n([nt])\to \ll(t)$ a.s.
, and consequently for $n$ large enough
$$\P(|\k_n([nt])-\ll(t)|>0)\leq\frac\eta 3.$$

On the other hand, thanks to Corollary \ref{coroll_alpha} we can find $u>0$ and $n_0\geq0$ such that for $n\geq n_0$, $\P(\ll^{-1}(u)<t)<\eta/3$, and $\frac{u}{\eps} \Big|\frac1{\alpha_n}-\frac1\alpha\Big|<\frac\eta 3$, where $\ll^{-1}(u):=\inf\{s\geq0,\ \ll(s)>u\}$ denotes the right inverse of $\ll$. Then for $n\geq n_0$,
\begin{align*}
 \P\left(\sum_{i=1}^{\ll(t)} |A_i-A_i^n|>\eps\right) & \leq \P\left(\sum_{i=1}^{u} |A_i-A_i^n|>\eps\right)+\P(\ll^{-1}(u)<t) \\
						    & \leq \frac{\E(\sum_1^u|A_i^n-A_i|)}{\eps}+\frac{\eta}{3} \\
						    & \leq \frac u\eps \Big|\frac1{\alpha_n}-\frac1\alpha\Big|+\frac\eta3\\
						    & \leq \frac{2\eta}{3},
\end{align*}
where the second inequality is obtained from an appeal to Markov's inequality. We conclude that $\underset{n\to\infty}\lim \P(|K_n[nt]-L_X(t)|>\eps)=0$.
\end{demo}

Next let us turn our attention back to our sequence of spectrally positive Lévy processes $(\tZ_n)$ converging to a Lévy process $Z$ with infinite variation. 

\begin{proposition} \label{CD_th3} %analogue théorème 3
If the convergence of $\tZ_n$ to $Z$ holds a.s., then for all $t\geq0$, we have convergence in probability of $L_n(t)$ towards $L(t)$.
\end{proposition}

\begin{demo}
 As in \cite{CD}, for all $n\geq0$, we consider the sequence of random walks $(\Snk)_{k\geq0}$ defined by $\Snk(j)=\tZ_n(j/k)$ for all $j\geq0$, so that as $k\to\infty$,
$$(\Snk([kt]))_{t\geq0}\to \tZ_n\ \text{a.s.}$$
As previously, each random walk $\Snk$ is endowed, independently of the others, with a Poisson process $N^{n,k}$ with parameter $\alpha_{n,k}:=\exp\{\sum_{i\geq1} \frac1i e^{-i/k}\P(\Snk(i)>0)\}$. We will use the obvious notation with subscript $n,k$ for all the quantities defined earlier involved by $\Snk$.\\
Fix $\eps>0$. From Proposition \ref{CD_th2}, we can find some sequence of integers $(k_n)_{n\geq1}$ such that, as $n\to\infty$,
$$(\Snkn([k_nt]))_{t\geq0}\to Z\ \text{a.s.}$$
and
$$\P(|K_{n,k_n}[k_nt]-L_n(t)|>\eps) \to 0.$$
We have 
\begin{align*}
 \P(| L_n(t)-L(t)|>3\eps) \leq & \ \ \P(|L_n(t)-K_{n,k_n}[k_nt]|>\eps) \\
			       & +\P(|K_{n,k_n}[k_nt]-\frac1{\alpha_{n,k_n}}\k_{n,k_n}([k_nt])|>\eps)\\
			       & +\P(|\frac1{\alpha_{n,k_n}}\k_{n,k_n}([k_nt])-L_t|>\eps) .
\end{align*}
We chose the subsequence $(k_n)$ such that the first term in the sum goes to $0$ as $n\to\infty$. The a.s. convergence of $(\Snkn)$ towards the Lévy process $Z$, for which the state $0$ is regular for $(0,\infty)$, allows us to apply Theorem 2 in \cite{CD} to get the convergence towards $0$ of the last term in the sum. \\
It remains to prove that $K_{n,k_n}[k_nt]-\frac1{\alpha_{n,k_n}}\k_{n,k_n}([k_nt])$ converges in probability to $0$ as $n\to\infty$. Recall that for all $n,j\geq0$, $\k_{n,k_n}(\t_{n,k_n}(j))=j$. Thus for all $j\geq0$, we can write 
\begin{align*}
\P(|K_{n,k_n}&[k_nt]  -\frac1{\alpha_{n,k_n}}\k_{n,k_n}([k_nt])|>\eps) \\
& = \P\left(\sum_{i=1}^{\k_{n,k_n}[k_nt]} \Big|A^{n,k_n}_i-\frac1{\alpha_{n,k_n}}\Big|>\eps\right)\\
& \leq \P\left(\sum_{i=1}^{[\alpha_{n,k_n}jt]} \Big|A^{n,k_n}_i-\frac1{\alpha_{n,k_n}}\Big|>\eps\right)+\P(\t_{n,k_n}[\alpha_{n,k_n}jt]<k_nt)\\
& \leq \frac{[\alpha_{n,k_n}jt]}{\eps^2 \alpha_{n,k_n}^2}+\P(\t_{n,k_n}[\alpha_{n,k_n}jt]<k_nt),
\end{align*}
the last inequality coming from the Bienaymé-Tchebitchev's inequality. From Remark 1 in \cite{CD}, we know that $\underset{n\to\infty}{\lim} \alpha_{n,k_n} =+\infty$. Thus letting first $n$ tend to $\infty$, we have that  $\frac{[\alpha_{n,k_n}jt]}{\eps^2 \alpha_{n,k_n}^2}$ goes to 0, and $\P(\t_{n,k_n}[\alpha_{n,k_n}jt]<k_nt)$ tends to $\P(\Linv(jt)<t)$ according to Theorem 1 in \cite{CD}. This last quantity now goes to $0$ as $j\to\infty$, and we completed the proof.
\end{demo}

\begin{corollaire}\label{CD_coroll2}
 The sequence $(\tZ_n,L_n,\Linv_n,\Han)$ converges as $n\to\infty$, in the sense of the finite dimensional distributions, to the process $(Z,L,\Linv,\Ha)$.
\end{corollaire}

\begin{demo} 
By Skorokhod's representation we may suppose that the convergence of $\tZ_n$ towards $Z$ holds a.s. Now Proposition \ref{prop_cv_H1} and Theorem \ref{th_cv_H} ensure the convergence in law of each coordinate, which provides the tightness of the quadruplet. Then, proving the a.s. convergence of the finite dimensional marginals will be sufficient to establish the corollary. \\
Now fix $t>0$. From Proposition \ref{CD_th3} we know that there exists some sequence of integers $k_n$, going to $\infty$ as $n\to\infty$, such that $L_{k_n}(t)$ tends to $L(t)$ a.s. From the definition of the inverse local time as a first passage time, and noting that $L$ has no fixed time of discontinuity, the latter convergence implies that of $\Linv_{k_n}(t)$ to $\Linv(t)$ a.s., by virtue of Proposition VI.2.11 in \cite{JS}. \\ 

As said in \cite{CD}, $\Linv(t)$ is an announceable stopping time (here $0$ is regular for $Z$ for $(0,\infty)$), so that from an appeal to Exercise 3 in \cite{B}, we get that $Z$ is a.s. continuous at time $\Linv(t)$. According to VI.2.3 in \cite{JS}, for all (possibly random) continuity point $u$ of $Z$, we have $\tZ_{k_n}(u)\to Z(u)$ a.s. as $n\to\infty$, and hence the sequence $(\tZ_{k_n}(t),L_{k_n}(t),\Linv_{k_n}(t),\tZ_{k_n}(\Linv_{k_n}(t)))$ converges a.s. as $n\to\infty$ towards $(Z(t),L(t),\Linv(t),\Ha(t))$.\\

Finally, taking any sequence of times $t_1<t_2<\ldots<t_j$, we can find a sequence $k'_n$ of integers tending to $\infty$ as $n\to\infty$, such that $\big((\tZ_{k'_n}(t_i),L_{k'_n}(t_i),\Linv_{k'_n}(t_i)),\Ha_{k'_n}(t_i),1\leq i\leq j\big)$ converges a.s. towards $\big((Z(t_i),L(t_i),\Linv(t_i),\Ha(t_i)),1\leq i\leq j\big)$ as $n\to\infty$, which ends the proof.
\end{demo}

\noindent\emph{Proof of Proposition \ref{prop_appendix} :} \\ 
According to Assumption \HypLevy\ (resp. Proposition \ref{CD_th3}, Proposition \ref{prop_cv_H1}, Remark \ref{remark_cv_H-}), we know that $\tZ_n$ (resp. $L_n$, $\Linv_n$, $\Han$, $\Hbn$) converges in distribution towards $Z$ (resp. $L$, $\Linv$, $\Ha$, $\Hb$). Therefore, each of these sequences is tight, and in consequence the sequence $(\tZ_n,L_n,\Linv_n,\Han,\Hbn)$ is tight. From Corollary \ref{CD_coroll2} we then get the joint convergence in distribution of $(\tZ_n,L_n,\Linv_n,\Han)$ towards $(Z,L,\Linv,\Ha)$, and moreover, the tightness ensures the existence of a subsequence $(k_n)$ such that $(\tZ_{k_n},L_{k_n},\Linv_{k_n},\Ha_{k_n},\Hb_{k_n})$ converges in distribution to $(\tilde Z,\tilde L,\tilde \Linv,\tHa, \tHb)$, with $(\tilde Z,\tilde L,\tLinv,\tHa)\stackrel{(d)}{=}(Z,L,\Linv,\Ha)$ and $\tHb\stackrel{(d)}{=}\Hb$. By virtue of the Skorokhod representation theorem, we can suppose that this convergence holds a.s. \\ 

The processes $\tHb$ and $\Hb$ are two subordinators and are equal in law, thus their continuous parts (which are deterministic drifts) are equal in law and therefore, almost surely. Consider now the jump part of $\tHb$. For all $\eps>0$ and $y\in\D(\R_+)$, define 

$$U(y,\eps):=\{u>0,\ |\Delta y(t)|=u \text{ for some } t\},$$
and
$$t^0(y,\eps):=0, \text{ and }  \forall p\geq0,\ t^{p+1}(y,\eps):=\inf\{t>t^p(y,\eps),\ |\Delta y(t)|>\eps\}.$$

Proposition VI.2.7 in \cite{JS} ensures for all $p\geq0$ that the mapping $y\mapsto t^p(y,\eps)$ (resp. $y\mapsto\Delta y(t^p(y,\eps))$) is continuous on $\D(\R_+)$ at each point $y$ such that $\eps\notin U(y,\eps)$ (resp. $\eps\notin U(y,\eps)$ and $t^p(y,\eps)<\infty$).  

Now we know that $t^p(Z,\eps)$, $t^p(\Ha,\eps)$, $t^p(\Hb,\eps)$ are finite a.s., and $\mu^+$, $\mu^-$ have no atoms. Moreover, since $\L$ is a $\sigma$-finite measure on $\R_+^*$, there exists a sequence $(\eps_m)_{m\geq1}$ of positive real numbers, which vanishes as $m\to\infty$, such that $\L(\{\eps_m\})=0$. As a consequence, for all $m\geq1$, the functions $y\mapsto t^p(y,\eps_m)$ and $y\mapsto\Delta y(t^p(y,\eps_m))$ are a.s. continuous w.r.t the distribution of $Z$, $\Ha$ and $\Hb$. Along with Proposition VI.2.1 of \cite{JS}, this gives for all $m\geq1$ the following almost sure convergence as $n\to\infty$:
 \begin{align*}
(t^p(\Ha_{k_n},\eps_m)&,\Delta \Ha_{k_n}(t^p(\Ha_{k_n},\eps_m)),\Delta \tZ_{k_n}(\Linv_{k_n}(t^p(\Ha_{k_n},\eps_m))),\Delta \Hb_{k_n}(t^p(\Ha_{k_n},\eps_m))) \\
&\xrightarrow{a.s.}(t^p(\tHa,\eps_m),\Delta \tHa(t^p(\tHa,\eps_m)),\Delta \tilde Z(\tLinv(t^p(\tHa,\eps_m))),\Delta \tHb(t^p(\tHa,\eps_m))).
 \end{align*}
Now with probability one the jumping times of $\Ha_n$ are exactly those of $\Hb_n$, and for all $t>0$, $\Delta \Hbn(t)=\Delta \tZ_n(\Linv_n(t))-\Delta \Han(t)$ a.s. Therefore letting now $m\to\infty$, we get :
\[\sum_{s\leq t} \Delta \tHb(s) = \sum_{s\leq t} (\Delta \tilde Z(\tLinv(s))-\Delta\tHa(s))\ \text{ a.s.}\]
As a consequence, we have $(\tilde Z,\tilde L,\tLinv,\tHa,\tHb)\stackrel{(d)}{=}(Z,L,\Linv,\Ha,\Hb)$, and then we get the convergence in distribution of $(\tZ_n,L_n,\Linv_n,\Han,\Hbn)$ towards $(Z,L,\Linv,\Ha,\Hb)$.
\hfill $\square$

\begin{remerciements}
 I would like to thank my supervisor, Amaury Lambert, for his very helpful advice and encouragement, and a reviewer for his careful reading of the paper.
\end{remerciements}

\end{document}